\pdfoutput=1
\RequirePackage{ifpdf}
\ifpdf 
\documentclass[pdftex]{sigma}
\else
\documentclass{sigma}
\fi

\numberwithin{equation}{section}
\numberwithin{theorem}{section}
\numberwithin{proposition}{section}
\numberwithin{lemma}{section}
\numberwithin{corollary}{section}
\numberwithin{conjecture}{section}
\numberwithin{definition}{section}
\numberwithin{example}{section}
\numberwithin{remark}{section}

\usepackage{amsxtra}
\usepackage{eucal}
\newcommand{\w}{\widetilde}

\newcommand{\cZ}{\mathcal Z}

\newcommand{\g}{{\mathfrak{g}}}

\newcommand{\Z}{{\mathbb Z}}

\newcommand{\N}{{\mathbb N}}

\newcommand{\bm}{{\mathbf m}}
\newcommand{\bn}{{\mathbf n}}
\newcommand{\bw}{{\mathbf w}}
\newcommand{\bx}{{\mathbf x}}
\newcommand{\by}{{\mathbf y}}
\newcommand{\bs}{{\mathbf s}}
\newcommand{\be}{{\mathbf e}}

\newcommand{\al}{{\alpha}}

\begin{document}

\allowdisplaybreaks

\renewcommand{\thefootnote}{$\star$}

\renewcommand{\PaperNumber}{014}

\FirstPageHeading

\ShortArticleName{$Q$-system Cluster Algebras, Paths and Total Positivity}

\ArticleName{$\boldsymbol{Q}$-system Cluster Algebras, Paths \\ and Total Positivity\footnote{This paper is a
contribution to the Proceedings of the Workshop ``Geometric Aspects of Discrete and Ultra-Discrete Integrable Systems'' (March 30 -- April 3, 2009, University of Glasgow, UK). The full collection is
available at
\href{http://www.emis.de/journals/SIGMA/GADUDIS2009.html}{http://www.emis.de/journals/SIGMA/GADUDIS2009.html}}}

\Author{Philippe DI FRANCESCO~$^\dag$ and Rinat KEDEM~$^\ddag$}

\AuthorNameForHeading{P. Di Francesco and R. Kedem}

\Address{$^\dag$~Institut de Physique Th\'eorique du Commissariat \`a l'Energie Atomique,
Unit\'e de Recherche\\
\hphantom{$^\dag$}~associ\'ee du CNRS, CEA Saclay/IPhT/Bat 774, F-91191 Gif sur Yvette Cedex,
France}

\EmailD{\href{mailto:philippe.di-francesco@cea.fr}{philippe.di-francesco@cea.fr}}
\URLaddressD{\url{http://ipht.cea.fr/en/Phocea/Pisp/visu.php?id=14}}

\Address{$^\ddag$~Department of Mathematics, University of Illinois
  Urbana, IL 61801, USA}
\EmailD{\href{mailto:rinat@illinois.edu}{rinat@illinois.edu}}
\URLaddressD{\url{http://www.math.uiuc.edu/~rinat/}}

\ArticleDates{Received October 15, 2009, in f\/inal form January 15, 2010;  Published online February 02, 2010}

\Abstract{In the f\/irst part of this paper, we provide a concise review of our
  method of solution of the $A_r$ $Q$-systems in terms of the
  partition function of paths on a weighted graph. In the second part,
  we show that it is possible to modify the graphs and transfer
  matrices so as to provide an explicit connection to the theory of
  planar networks introduced in the context of totally positive
  matrices by Fomin and Zelevinsky. As an illustration of the further
  generality of our method, we apply it to give a simple
  solution for the rank~2 af\/f\/ine cluster algebras studied by Caldero
  and Zelevinsky.}


\Classification{05E10; 13F16; 82B20}

\setcounter{tocdepth}{2}\tableofcontents

\section{Introduction}

Discrete dynamical systems may take the form of recursion relations
over a discrete time variable, describing the evolution of relevant
physical quantities. Within this framework, of particular interest are
the discrete integrable recursive systems, for which there exist
suf\/f\/iciently many conservation laws or integrals of motion, so that
their solutions can be expressed in terms of some initial data.
Interesting examples of such systems of non-linear integrable
recursion relations arise from matrix models used to generate random
surfaces, in the form of discrete Toda-type equations
\cite{AVM,KKN,KLWZ}.  More recently, a combinatorial study of
intrinsic geometry in random surfaces has also yielded a variety of
integrable recursion relations, also related to discrete spatial
branching processes \cite{GEOD}.

We claim that other fundamental examples are provided by the so-called
$Q$-systems for Lie groups, introduced by Kirillov and Reshetikhin
\cite{KR} as combinatorial tools for addressing the question of
completeness of the Bethe ansatz states in the diagonalization of the
Heisenberg spin chain based on an arbitrary Lie algebra. We proved
integrability for these systems in the case of $A_r$ in
\cite{DFK08a}. In the case of other Dynkin diagrams, evidence suggests
integrability still holds.

The $Q$-system with special (singular) initial conditions was originally
introduced \cite{KR} as the recursion relation satisf\/ied by the characters of special
f\/inite-dimensional modules of the Yangian~$Y(\g)$, the so-called Kirillov--Reshetikhin modules.
Remarkably, in the case $\g=A_r$, the same recursion relation also appears
in other contexts, such as Toda f\/lows in Poisson geometry~\cite{GSV}, preprojective algebras~\cite{GLS} and canonical bases~\cite{BER}.

In \cite{DFK08a}, we used methods from statistical mechanics to study
the solutions of the $Q$-system associated with the Lie algebra~$A_r$,
for f\/ixed but arbitrary initial conditions.  Our approach starts with
the explicit construction of the conserved quantities of the system,
which appear as coef\/f\/icients in a linear recursion relation satisf\/ied
by $Q$-system solutions.  These are f\/inally used to reformulate the
solutions in terms of partition functions for weighted paths on
graphs, the weights being entirely expressed in terms of the initial
data.

\looseness=-1
Note that there is a choice of various sets of $2r$ variables which
constitute an initial condition
f\/ixing the solutions of the $Q$-system recursion relation, a choice
parametrized by Motzkin paths of length $r$. This set of initial
variables determines the graphs and weights which solve the
problem. This is the key point addressed in~\cite{DFK08a}, and is
related to the formulation as a cluster algebra.

Cluster algebras \cite{FZI} are another form of discrete dynamical
systems. They describe a specif\/ic type of evolution, called mutation,
of a set of variables or cluster seed. Mutations are rational,
subtraction-free expressions.  This type of structure has
proved to be very universal, and arises in many dif\/ferent mathematical
contexts, such as total positivity \cite{FZposit,FZpositII}, quiver
categories~\cite{BK}, Teichm\"uller space geometry~\cite{TRIANG},
Somos-type sequences~\cite{FZLaurent}, etc.

Cluster algebras have the property that any cluster variable is
expressible as a Laurent polynomial of the variables in any other
cluster in the algebra. It is conjectured that these Laurent
polynomials have nonnegative coef\/f\/icients \cite{FZI} (the positivity
conjecture).  This property has only been proved in a few
context-specif\/ic cases so far, such as f\/inite type acyclic case \cite{FZIV},
af\/f\/ine type acyclic case~\cite{ARS}, or clusters arising from surfaces~\cite{MSW}.
The $Q$-system solutions for $A_r$ are also known to form a subset of the cluster
variables in the cluster algebra introduced in~\cite{Ke07} (a result
later generalized to all simple Lie algebras in \cite{DFK08}). In~\cite{DFK08a}, we interpreted the solutions of the~$A_r$ $Q$-system in terms of
partition functions of paths on graphs with positive weights: this proved positivity
for the corresponding subset of clusters.  Moreover, we
obtained explicit expressions, in the form of f\/inite continued
fractions, for these cluster variables. We review our results in the
f\/irst part of this paper.

Of particular interest to us is the connection
of cluster algebras to total positivity. Fomin and Zelevinsky
\cite{FZposit} expressed a parametrization of totally positive
matrices in terms of electrical networks, and established total
positivity criteria based on relations between matrix minors,
organized into a cluster algebra structure. In this paper we show the
explicit connection of their construction to the $Q$-system solutions.

\looseness=1
In this paper, we review the methods and results of \cite{DFK08a} in a
more compact and hopefully accessible form. We f\/irst apply this method
to the case of rank 2 af\/f\/ine cluster algebras, studied by Caldero and
Zelevinsky~\cite{CZ}. These are the cluster algebras which arise from
the Cartan matrices of af\/f\/ine Kac--Moody algebras. We obtain a simple
explicit solution for the cluster variables in terms of initial data.
We then proceed to describe the general solution of the $A_r$ $Q$-system
using the same methods. In particular, we obtain an explicit formula for
the fundamental cluster variables, which generalizes
the earlier results of~\cite{CZ} to higher rank.
Finally, we make
the explicit connection between the path interpretation of
the solutions of the $Q$-system
and a~subclass of the totally positive matrices of~\cite{FZposit} and their associated
electrical networks.

More precisely, it turns out that the generating function for the
family of cluster variables~$R_{1,n}$ $(n\in \N)$ of the $A_r$
$Q$-system is given by the resolvent of the transfer matrix $T_\bm$
associated with a graph $\Gamma_\bm$ for some seed variables
associated with the Motzkin path $\bm$. We show that it is possible to
locally modify the graph without changing the path generating
function, so that we obtain a transfer matrix of smaller size $r+1$,
equal to the rank of the algebra. From there, there is a
straightforward identif\/ication with the networks associated with
totally positive matrices of special type, related to the Coxeter
double Bruhat cells of \cite{GSV}.

The advantage of our approach is that we have explicit expressions
for the cluster variables in terms of any mutated cluster seed
parametrized by a Motzkin path $\bm$.

The paper is organized as follows. In Section~\ref{hard}, we explain the basic
tools from statistical mechanics which we use, the partition functions
of hard particles on graphs and the equivalent partition functions of
paths on weighted dual graphs.

For illustration, in Section~\ref{section3}, we use
this to give the explicit expression for the gene\-rating function of
the cluster variables of rank~2 cluster algebras corresponding to
af\/f\/ine Dynkin diagrams, a problem extensively studied in \cite{SZ,CZ,PM}.

In Section~\ref{qsysappli}, we review our solution~\cite{DFK08a} of the $A_r$
$Q$-system for the simplest choice of initial variables. This solution
uses the partition functions of Section~\ref{hard} as well as the theorem of
\cite{LGV1,LGV2} relating the partition function of non-intersecting
paths to determinants of partition functions of paths. We also
introduce a new notion in this section, the ``compactif\/ication'' of
the graph, which gives a new transfer matrix which is associated with
the same partition function. This is a key tool in making the
connection with totally positive matrices.

Section~\ref{generalpaths} generalizes the results of Section~\ref{qsysappli}, and we give
expressions for the $Q$-system solutions in terms of any set of
cluster variables in a fundamental domain. We also introduce the
notion of ``strongly non-intersecting paths'' and a generalization of
\cite{LGV1,LGV2}.  This gives generating functions for the cluster
variables in terms of the other seeds in the cluster algebra, and also
provides a proof of the positivity conjecture \cite{FZI}. This section
is a quick review of the results of \cite{DFK08a}.

In Section~\ref{section6}, we extend the graph compactif\/ication procedure of
Section~\ref{qsysappli} to the other graphs, corresponding to mutated cluster
variables, introduced in Section~\ref{generalpaths}. This yields transfer matrices of
size $(r+1)\times (r+1)$, the resolvent of which is an alternative
expression for the generating function for cluster variables. Finally,
in Section~\ref{totalpos}, we use this to give the explicit relation to totally
positive matrices \cite{FZposit,GSV}.

\section{Partition functions}\label{hard}

\subsection[Hard particles on $G_r$]{Hard particles on $\boldsymbol{G_r}$}

\subsubsection{Vertex-weighted graphs} Let $G$ be a f\/inite graph with
$N$ vertices
labeled $1,2,\dots ,N$, and single, non-oriented edges connecting some
vertices.  The adjacency matrix $A^G$ of the graph $G$ is the
$N\times N$ matrix with entries $A^G_{i,j}=1$ if vertex $i$ is
connected by an edge to a vertex~$j$, and $0$ otherwise.
To each vertex~$i$ is associated a positive
weight $y_i$.

\subsubsection{Conf\/igurations}

A conf\/iguration $C$ of $m$ hard particles on $G$ is a subset of
$\{1,\dots ,N\}$ containing $m$ elements, such that $A^G_{i,j}=0$ for all
$i,j\in C$.

This is called a hard particle conf\/iguration, because if we view the
elements of $C$ to be the vertices occupied by particles on the graph, the
condition $A^G_{i,j}=0$ for $i,j\in C$ enforces the rule that two
neighboring sites cannot be occupied at the same time. Each vertex
can be occupied by at most one particle.

We denote by $\mathcal C^G_m$ the set of all hard particle
conf\/igurations on $G$ with $m$ particles.

\subsubsection{Partition function}

The weight of a conf\/iguration $C$ is the product of all the weights
associated with the elements of $C$. That is,
\[
w_C = \prod_{i\in C} y_i.
\]

The partition function $Z^{G}_m$ for $m$ hard particles on $G$ is
the sum over all conf\/igurations $\mathcal C^G_m$ of the corresponding
weights:
\begin{equation*}
Z^{G}_m(\by)=\sum_{C\in \mathcal C_m^{G}} w_C.
\end{equation*}

\begin{figure}[t]
\centerline{\includegraphics{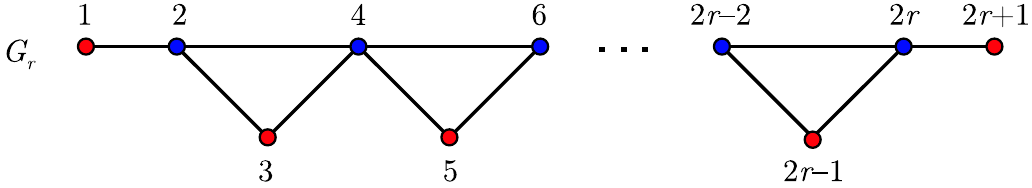}}
\caption{The graph $G_r$, with $2r+1$ vertices labeled $i=1,2,\dots ,2r+1$.}
\label{fig:graphgr}
\end{figure}

\subsubsection[The graph $G_r$]{The graph $\boldsymbol{G_r}$}

A basic example is the graph $G=G_r$ of Fig.~\ref{fig:graphgr}. It
has $2r+1$
vertices and $3r-1$ edges, with the (symmetric) adjacency matrix def\/ined by
\begin{gather*}
 A_{1,2}=A_{2r,2r+1}=A_{2 i, 2i+1}=A_{2i, 2i+2}=A_{2i+1,2i+2}=1,\qquad i\in
\{1,\dots ,r-1\}.
\end{gather*}

When $r=0$, $G_0$ is reduced to a single vertex labeled~$1$, while for
$r=1$, $G_1$ is a chain of three vertices $1$, $2$, $3$ with the two edges~$(1,2)$ and $(2,3)$.

\begin{example} \label{zeroonehp}
The non-vanishing partition functions $Z_m^{G_r}$ for hard particles on~$G_r$,
$r=0,1$
are
\begin{gather}\label{initZ}
Z_0^{G_0}=1,\qquad Z_1^{G_0}=y_1, \qquad Z_0^{G_1}
=1,\qquad Z_1^{G_1}=y_1+y_2+y_3, \qquad  Z_2^{G_1}=y_1 y_3.
\end{gather}
\end{example}

\subsubsection{Recursion relations for the partition function}
The transfer matrices $Z_m^{G_r}$ satisfy recursion relations in the
index $r$. They are obtained by considering the possible occupancies
of the vertices $2r+1$ and $2r$:
\begin{equation}\label{recuHO}
Z_m^{G_r}=Z_m^{G_{r-1}}+ y_{2r+1} Z_{m-1}^{G_{r-1}}
+y_{2r} Z_{m-1}^{G_{r-2}}, \qquad r\geq 2,  \quad m\geq 0.
\end{equation}
For example, $Z_0^{G_r}=1$ is the partition function of the empty conf\/iguration
and $Z_{r+1}^{G_r}=\prod\limits_{i=0}^r y_{2i+1}$ for the maximally occupied conf\/iguration.

\begin{figure}[t]
\centerline{\includegraphics{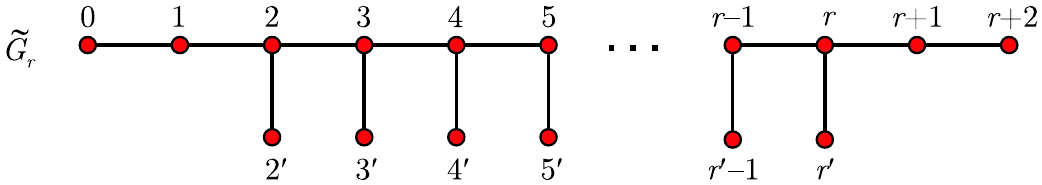}}
\caption{The graph ${\tilde G_r}$, dual to $G_r$,
with  total of $2r+2$ vertices.}
\label{fig:dualgr}
\end{figure}

\subsection{Transfer matrix on the dual graph}
Associated to the graph $G_r$, there is a dual graph ${\widetilde
  G}_r$, as in Fig.~\ref{fig:dualgr}.

It is dual in the sense that $G_r$ is the medial graph of ${\widetilde
  G}_r$: each edge
of ${\widetilde G}_r$ corresponds to a~vertex of $G_r$, and any two edges
of ${\widetilde G}_r$ share a vertex if\/f the corresponding vertices of
$G_r$ are adjacent.

We f\/ix the labeling so that the correspondence is between edges of
${\widetilde G}_r$ and vertices of $G_r$ is:
\begin{enumerate}\itemsep=0pt
\item[1)] edge $(k,k')$ of ${\widetilde G}_r$ corresponds to the vertex
  $2k-1$ of $G_r$, where $k=2,3,\dots ,r$;
\item[2)] edge $(k,k+1)$ of ${\widetilde G}_r$ corresponds to the vertex
  $2k$ of $G_r$, where $k=1,2,\dots ,r$;
\item[3)] edge $(0,1)$ of ${\widetilde G}_r$ corresponds to the vertex $1$
  of $G_r$;
\item[4)] edge $(r+1,r+2)$ of ${\widetilde G}_r$ corresponds to the vertex $2r+1$
  of $G_r$.
\end{enumerate}

\subsubsection[Transfer matrix on $\widetilde G_r$]{Transfer matrix on $\boldsymbol{\widetilde G_r}$}
\begin{figure}[t]
\centerline{\includegraphics{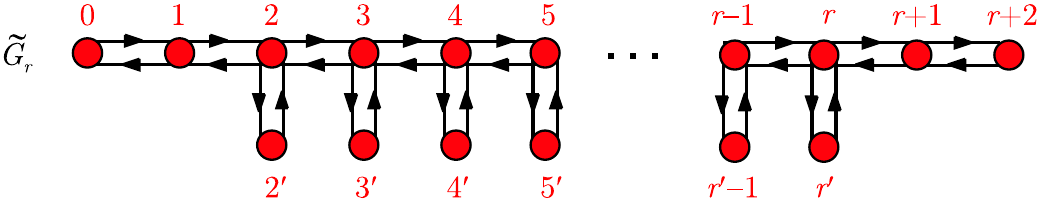}}
\caption{The graph ${\tilde G_r}$ with oriented edges drawn in.}
\label{fig:dualgrter}
\end{figure}
We choose an ordering of the $2r+2$ vertices of ${\widetilde G}_r$ to be
$0<1<2<2'<3<3'<\dots <r<r'<r+1<r+2$, and consider this to be an index set. We
construct the $(2r+2)\times (2r+2)$
transfer matrix $T_r(t\by)$ with these indices. Its entries are $0$ except for:
\begin{gather*}
 (T_r)_{k',k}=1,\quad (T_r)_{k,k'}=ty_{2k-1},\qquad
k=2,3,\dots ,r;\\
 (T_r)_{k+1,k}=1,\qquad (T_r)_{k,k+1}=ty_{2k},\qquad  k=1,2,\dots ,r;\\
  (T_r)_{1,0}=(T_r)_{r+2,r+1}=1,\qquad (T_r)_{0,1}=ty_1,\qquad
(T_r)_{r+1,r+2}=ty_{2r+1}.
\end{gather*}
In matrix form,
\begin{equation}\label{transmat}
T_r(t\by)=
 \left(
\begin{matrix}
0 & ty_1 & 0     & 0    & 0    & \cdots & 0 & 0 & 0 & 0 \\
1 & 0     & ty_2 & 0    & 0    & \cdots & 0 & 0 & 0 & 0 \\
0 & 1     & 0    &  ty_3&ty_4 &  \cdots & 0 & 0 & 0 & 0 \\
0 & 0     & 1    &  0    &  0   & \cdots & 0 & 0 & 0 & 0 \\
0 & 0     & 1    &  0    &  0   & \cdots & 0 & 0 & 0 & 0 \\
\vdots &\vdots     &     &     &     & \ddots &  &  &  & \vdots \\
0 & 0     & 0    &  0    &  0   & \cdots & 0 & ty_{2r-1} & ty_{2r} & 0 \\
0 & 0     & 0    &  0    &  0   & \cdots & 1 & 0           & 0         & 0 \\
0 & 0     & 0    &  0    &  0   & \cdots & 1 & 0           & 0         & ty_{2r+1} \\
0 & 0     & 0    &  0    &  0   & \cdots & 0 & 0           & 1         & 0
\end{matrix}
\right).
\end{equation}
This matrix is a weighted adjacency matrix for the graph $\widetilde G_r$ as drawn in Fig.~\ref{fig:dualgrter}. Each edge of $\widetilde G_r$
corresponds to two oriented edges pointing in opposite directions,
with the weights in the transfer matrix corresponding to those oriented
edges. The element of the transfer matrix indexed by $i$, $j$ is the
weight corresponding to the edge $j\to i$.

\begin{lemma}
The generating function for hard particle partition functions on $G_r$
with weights $t y_i$ per particle at vertex $i$ is
\begin{equation}\label{hardobT}
Z^{G_r}(t\by):=\sum_{m=0}^{r+1} t^mZ_m^{G_r}=\det(I-T_r(-t \by))
\end{equation}
with $T_r$ as in equation~\eqref{transmat}.
\end{lemma}

\begin{proof}
Expanding the determinant
$D_r(t)=\det(I-T_r(-t\by))$ along the last column, one obtains
the recursion relation:
\begin{equation*}
D_r(t)= D_{r-1}(t)+ty_{2r+1}D_{r-1}(t)+ty_{2r}D_{r-2}(t),
\qquad r\geq 2.
\end{equation*}
Let $D_{r,m}$ denote
the coef\/f\/icient of $t^m$ in $D_r(t)$. It satisf\/ies the same
recursion relation as equation~\eqref{recuHO} for $Z_m^{G_r}$, with
$D_0(t)=1+ty_1$
and $D_1(t)=1+t(y_1+y_2+y_3)+t^2y_1y_3$, in agreement with the initial
values of $Z^{G_r}(t\by)$, $r=0,1$ in equation \eqref{initZ}. Thus,
$D_{r,m}=Z_m^{G_r}$ for all $r$, $m$.
\end{proof}

\begin{example}\label{examgone}
For the case $r=1$, $G_1$ is a chain of 3 vertices $1$, $2$, $3$. The dual ${\widetilde G}_1$
is a chain of three edges connecting four vertices $0$, $1$, $2$, $3$, and the transfer matrix
is:
\begin{equation*}
T(t\by)=\left( \begin{matrix}
0 & ty_1 & 0 & 0 \\
1 & 0 & ty_2 & 0 \\
0 & 1 & 0 & ty_3 \\
0 & 0 & 1 & 0
\end{matrix}\right).
\end{equation*}
One checks that $\det(I-T(-t\by))=1+t(y_1+y_2+y_3)+t^2 y_1y_3$, in agreement
with the partition functions $Z_i^{G_1}$ of equation~\eqref{initZ}.
\end{example}

\subsection{Hard particles and paths}\label{hardtopath}
Let $G$ be a graph with oriented edges. For example, the graph
$\widetilde{G}_r$ of the previous section can be made into an oriented
graph, by taking each edge to be a
doubly-oriented edge. Then each oriented edge from $i$ to $j$ receives
a weight $w(i,j)$, which is the corresponding entry of the transfer
matrix $T(t\by)$ of equation \eqref{transmat}.
This may be interpreted as the
transfer matrix for paths on ${\widetilde G}_r$ as follows.

\subsubsection{The partition function of paths}\label{pathinv}

A path of length $n$ on an oriented graph $G$ is a sequence of
vertices of $G$, $P=(v_0,v_1,\dots ,v_n)$, such that there exists an arrow
from $v_i$ to $v_{i+1}$ for each $i$.  Let $\mathcal P_{a,b}^G(n)$ be the set
of all distinct paths of length $n$, starting at vertex $a$ and ending
at vertex $b$ on the graph $G$.

A path on a graph with weighted edges has a total weight which is the
product of the weights $w(v_i,v_{i+1})$ associated with
 the edges traversed in the path.

The
partition function ${\mathcal Z}_{a,b}^{G}(n)$ for $\mathcal
P_{a,b}^{G}(n)$ is
\begin{gather*}
{\mathcal Z}_{a,b}^{G}(n)=\sum_{P\in \mathcal P_{a,b}^G(n)} \prod_{i=0}^{n-1}
w(v_i,v_{i+1}).
\end{gather*}
Assuming $G$ is f\/inite, and labeling its vertices $i=1,2,\dots ,N$, let
us introduce the $N\times N$ transfer matrix $T$ with entries
$T_{i,j}=w(j,i)$, we have the following simple expression for
${\mathcal Z}_{a,b}^{G}(n)$:
\begin{equation*}
{\mathcal Z}_{a,b}^{G}(n)= \big( T^n \big)_{b,a}.
\end{equation*}

The matrix $T(t\by)$ of~\eqref{transmat} is then
the transfer matrix for paths on ${\widetilde G}_r$ with weights $1$
for edges pointing away from the vertex $0$, and $ty_i$ for the
$i$-th edge pointing to the origin. The partition function for
weighted paths of arbitrary length on ${\widetilde G}_r$ from the
vertex $0$ to itself is
\begin{equation}\label{genpath}
{\mathcal Z}^{{\widetilde G}_r}(t\by)
=\sum_{n\geq 0} \big( T(t\by)^n \big)_{0,0} =\big((I -T(t\by))^{-1}\big)_{0,0}.
\end{equation}

\begin{lemma}
The partition function of paths of arbitrary length from $0$ to $0$ on
the graph $\widetilde{G}_r$ is equal to
\begin{equation}\label{heaps}
{\mathcal Z}^{{\widetilde G}_r}(t\by)={Z^{G_r}(0,-ty_2,\dots ,-ty_{2r+1})\over
Z^{G_r}(-ty_1,-ty_2,\dots ,-ty_{2r+1})}.
\end{equation}
\end{lemma}
\begin{proof}
By def\/inition, \eqref{genpath} is the ratio of the $(0,0)$-minor of the matrix
$I-T(t\by)$ to its determinant.  Using equation \eqref{hardobT} and the
explicit form \eqref{transmat}, this immediately yields the relation.
\end{proof}

This relation between hard-particle partition functions and path
partition functions may be interpreted as a boson-fermion
correspondence, and is a particular case of Viennot's theory of heaps
of pieces \cite{HEAPS1,HEAPS2}.

\subsubsection{The path partition function as a continued
  fraction}\label{contifrac}

A direct way to compute ${\mathcal Z}^{{\widetilde G}_r}(t\by)$ in
equation \eqref{genpath} is by using Gaussian elimination on
$I-T(t\by)$ to bring it to lower-triangular form.  The resulting pivot
in the f\/irst row is $1/\big((I-T(t\by))^{-1}\big)_{0,0}$. If we do
this systematically, by left-multiplication by upper-triangular
elementary matrices, the result is
\begin{lemma}
\begin{equation}\label{continued}
{\mathcal Z}^{{\widetilde G}_r}(t\by)={1\over 1-t{y_1\over
    1-t{y_2\over 1-t y_3 -t{y_4\over
1-ty_5-t{y_6\over {\ddots \over 1-ty_{2r-1}-t{y_{2r}\over 1-ty_{2r+1}}}}}}}}.
\end{equation}
\end{lemma}

\subsubsection{Non-Intersecting paths: the LGV formula}

\begin{figure}[t]
\centerline{\includegraphics{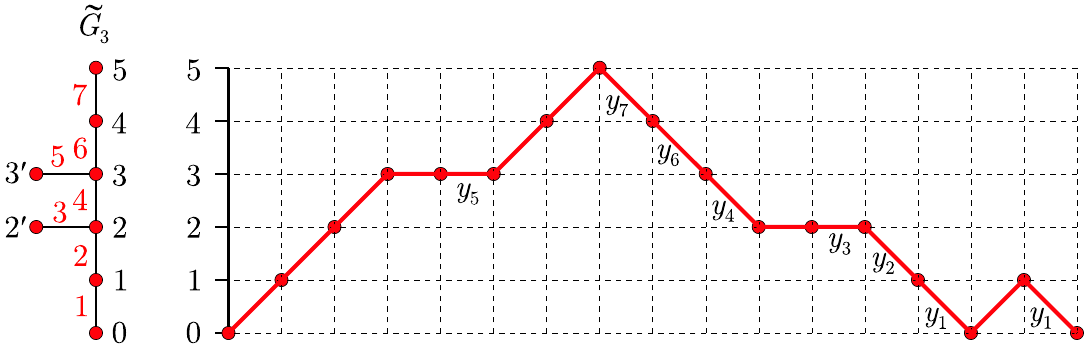}}
\caption{\small A path on ${\widetilde G}_3$ with $16$ steps.
The weights $y_i$ are associated to the descending steps $(1,-1)$ and to
the second half of the horizontal steps $(1,0)+(1,0)$ of the path, $i$ being
the label of the corresponding edge of ${\widetilde G}_3$ (see on left).
Here, the path receives the weight $y_1^2y_2y_3y_4y_5y_6y_7$.}\label{fig:grpath}
\end{figure}

We may represent paths of length $2n$ from vertex $0$ to itself on
${\widetilde G}_r$ as paths on the lattice $\Z_{\geq 0}^2$ (see
Fig.~\ref{fig:grpath} for an illustration). Such paths start at
$(0,0)$ and end at $(2n,0)$, and have the following possible steps:{\samepage
\begin{enumerate}\itemsep=0pt
\item[1)] to the northeast, $(j,k)\to(j+1,k+1)$, corresponding to the
  $j$th step in the path going from vertex $k$ to vertex $k+1$;
\item[2)] to the southeast, $(j,k)\to(j+1,k-1)$, corresponding to the
  $j$th step in the path going from vertex $k$ to vertex $k-1$;
\item[3)] to the east, $(j,k)\to(j+2,k)$, corresponding to the steps $k\to
k'\to k$ if $k\in\{2,3,\dots ,r\}$.
\end{enumerate}
Paths from the origin to the origin have an even number of steps, by
parity.}

We will need to consider the partition function of families of $\al$
non-intersecting paths on $\widetilde{G}_r$, ${\mathcal
  Z}_{\bs,\be}^{{\widetilde G}_r}$. Here, the f\/ixed starting points
are parametrized by $\bs=(s_1,\dots ,s_\al)$ and the endpoints by
$\be=(e_1,\dots,e_\al)$.

By non-intersecting paths, we mean paths that do not share any vertex.
We have the celebrated Lindstr\"om--Gessel--Viennot formula~\cite{LGV1,LGV2}
\begin{equation}\label{lgv}
{\mathcal Z}_{\bs,\be}^{{\widetilde G}_r}=\det_{1\leq i,j \leq \al} \,
{\mathcal Z}_{s_i,e_j}^{{\widetilde G}_r}.
\end{equation}
The determinant has the ef\/fect of subtracting the contributions from
paths that do intersect. This formula can be proved by expanding the
determinant as:
\begin{equation}\label{expandet}
\sum_{\sigma \in S_\al} {\rm sgn}(\sigma)\, \prod_{i=1}^\al
{\mathcal Z}_{s_{\sigma(i)},e_i}^{{\widetilde G}_r}.
\end{equation}
One then considers the involution $\varphi$
on families of paths, def\/ined by interchanging the beginnings of the two f\/irst paths that share
a vertex until the vertex, or by the identity if no two paths in the family intersect. It is clear
that when $\varphi$ does not act as the identity it relates two path conf\/igurations with opposite
weights in the expansion of the determinant, as the two starting points are switched by~$\varphi$,
hence these cancel out of the expansion~\eqref{expandet}.
We are thus left only with non-intersecting families,
all corresponding to $\sigma={\rm Id}$, hence all with positive weights.

\section{Application to rank 2 cluster algebras of af\/f\/ine type}\label{section3}

\subsection{Rank two cluster algebras}
In this section, we use the partition functions introduced in the
previous section to the problem of computing the cluster variables of
rank two cluster algebras of af\/f\/ine type with trivial coef\/f\/icients. This allows us to give
explicit, manifestly positive formulas for the variables, proving the
positivity of the variables in these cases.

\subsubsection{The recursion relations}

The rank 2 cluster algebras of af\/f\/ine type \cite{FZI,SZ} with trivial coef\/f\/icients may be reduced
to the following recursion relations for $n\in \Z$:
\begin{equation}\label{recuaff}
  x_{n+1}=\left\{ \begin{matrix} \dfrac{1}{x_{n-1}}(1+x_n^b) & {\rm if}\ \ n \ \ {\rm is}\ \ {\rm odd},\vspace{2mm}\\
      \dfrac{1}{x_{n-1}}(1+x_n^c) & {\rm if}\ \ n \ \ {\rm is}\ \ {\rm even}, \end{matrix}\right.
\end{equation}
where $b$, $c$ are two positive integers with $bc=4$, hence
$(b,c)=(2,2)$, $(4,1)$ or $(1,4)$. The aim is to
f\/ind an expression for $x_n$, $n\in \Z$ in terms of some initial data,
e.g.~$(x_0,x_1)$.

The connection to rank 2 af\/f\/ine Lie algebras is via the Cartan matrix
$\left(
  \begin{matrix}2 & -b\\ -c & 2\end{matrix} \right)$.

The cases $(4,1)$ and $(1,4)$ are almost equivalent: If $x_n(x_0,x_1)$
is a solution of the $(b,c)$ equation, then $x_{1-n}(x_1,x_0)$ is the
solution of the $(c,b)$ equation.  We may thus restrict ourselves to
the $(1,4)$ case.

However, the symmetry $x_n\leftrightarrow x_{1-n}$
changes the parity of $n$. Therefore we need to also
consider the dependence of $x_n$ in the ``odd'' initial data
$(x_1,x_2)$.

It turns out that the recursion relations \eqref{recuaff} are all integrable
evolutions. This allows us to compute the generating function for
$x_n$, $n\geq 0$. The result is a manifestly positive (f\/inite)
continued fraction.  In the light of the results of the previous
section, this allows to reinterpret $x_n$ as the partition function
for weighted paths on certain graphs. This path formulation gives yet
another direct combinatorial interpretation for the expression of
$x_n$ as a positive Laurent polynomial of the initial data, to be compared with
the approach of \cite{SZ,CZ} using quiver representations and that of~\cite{PM} using matchings of dif\/ferent kinds of graphs.

\subsection[The $(2,2)$ case: solution and path interpretation]{The $\boldsymbol{(2,2)}$ case: solution and path interpretation}\label{twotwo}

Consider the recursion relation
\begin{equation}\label{qsysaone}
x_{n+1}x_{n-1}=x_n^2+1
\end{equation}
with $\bx_0=(x_0,x_1)=(x,y)$. This is the $A_1$ case of the
renormalized $Q$-system considered in \cite{DFK08}.

Due to the symmetry
$n\leftrightarrow 1-n$ of the equation, the solution satisf\/ies
\[
x_n(x_0,x_1)=x_{1-n}(x_1,x_0)
\]
so we may restrict our attention to computing $x_n$, for $n\geq 0$.

\subsubsection{Constants of the motion}
Equation \eqref{qsysaone} is integrable. To see this, we rewrite
\eqref{qsysaone} as
\[
\varphi_n:=\left\vert \begin{matrix} x_{n-1} & x_n\\ x_n & x_{n+1}\\
\end{matrix}\right\vert=1.
\]
Then
\[
0=\varphi_n-\varphi_{n+1}=\left\vert \begin{matrix} x_{n-1}+x_{n+1} & x_n\\
    x_n+x_{n+2} & x_{n+1}\\ \end{matrix}\right\vert.
\]
We conclude
that there exists a constant $c$, independent of $n$, such that
$x_{n-1}+x_{n+1}=c x_n$ for all $n$.  Using equation \eqref{qsysaone},
\begin{equation}\label{iom}
c={x_{n-1}+x_{n+1}\over x_n}={x_{n+1}\over x_n}+{1\over x_nx_{n+1}}+{x_n\over x_{n+1}}
={x_1\over x_0}+{1\over x_0x_{1}}+{x_0\over x_{1}}.
\end{equation}
We interpret $c$ as an integral of motion of the three-term relation
\eqref{qsysaone}: all solutions of the latter indeed satisfy the two-term recursion
relation \eqref{iom}, for some ``integration constant''~$c$ f\/ixed by the initial data.

Note that $c$ coincides with the partition function
$Z_1^{G_1}=y_1+y_2+y_3$ for one hard particle on the graph $G_1$, with
weights
\begin{equation}\label{yweights}
y_1={x_1\over x_0},\qquad y_2={1\over x_0x_1},\qquad y_3={x_0\over
  x_1}.
\end{equation}
Note also that the only other non-vanishing hard particle
partition functions on $G_1$ are $Z_0^{G_1}=1$ and $Z_2^{G_1}=y_1y_3=1$.

\subsubsection[Generating function for $x_n$]{Generating function for $\boldsymbol{x_n}$}

Let $X(t)=\sum\limits_{n=0}^\infty t^n x_n$ be the generating function for
the variables $x_n$ with $n\geq 0$. Using
$x_{n+1}-cx_n+x_{n-1}=0$, we have by direct calculation
\begin{equation}\label{aoneres}
X(t)={x_0- (cx_0-x_1)t \over 1-c t+t^2}={x_0\over 1-t {y_1\over 1-t {y_2\over 1-t y_3}}},
\end{equation}
with $y_i$ as in \eqref{yweights}.
This gives $x_n$ as a manifestly positive Laurent polynomial of $(x_0,
x_1)$.  In fact, expanding the r.h.s.\ of~\eqref{aoneres}, we get
\begin{equation*}
X(t)=\sum_{p,q,\ell\geq 0}{p+q-1\choose q}{q+\ell-1\choose \ell}
t^{p+q+\ell}x_0^{1+\ell-p-q}x_1^{p-q-\ell},
\end{equation*}
from which $x_n$ is obtained by extracting the coef\/f\/icient of $t^n$.
This agrees with the result of \cite{CZ}.

It follows that the dependence of the variables $x_n$ on any other pair
$(x_k,x_{k+1})$ is also as a~positive Laurent polynomial. This is
clear from the translational invariance of the system:
\[
x_n(x_k,x_{k+1}) = x_{n-k}(x_0,x_{1})\big|_{x_0\mapsto x_k,x_1\mapsto x_{k+1}}.
\]

\subsubsection{Relation to the partition function of paths} Upon
comparing the continued fraction expression \eqref{aoneres} and
equation \eqref{continued}, we see that there is a path interpretation
to the variables $x_n$ as follows.

The denominator of the fraction
$X(t)$ is equal to the partition function $Z^{G_1}(-t \by)$
for hard particles on the graph $G_1$ introduced in Section~\ref{hard}, with the weights $-t y_i$ per particle at vertex~$i$.
Therefore,
\begin{equation*}
X(t)=x_0 \, {Z^{G_1}(0,-ty_2,-ty_3)\over Z^{G_1}(-t y_1,-ty_2,-ty_3)},
\end{equation*}
so $x_0^{-1}X(t)$ is equal the partition function for paths on the
graph ${\widetilde G}_1$ of Example \ref{examgone}, beginning and
ending at the vertex 0. The weights are as follows: The weight is
equal to 1 for any step away from vertex 0, and
and is equal to $t y_i$ per step from vertex $i$ to $i-1$.

\subsection[The $(1,4)$ case: solution and relation to paths]{The $\boldsymbol{(1,4)}$ case: solution and relation to paths}

We now consider the system:
\begin{gather*}
x_{2n}={1+x_{2n-1}\over x_{2n-2}},\nonumber \\
x_{2n+1}={1+x_{2n}^4\over x_{2n-1}}, \qquad n\in \Z.
\end{gather*}
We determine the dependence of the variables $x_n$ on two dif\/ferent
types of initial conditions:
\begin{gather}
{\rm case}\ 0:\ (x_0,x_1), \nonumber \\
{\rm case}\ 1: \ (x_1,x_2).\label{cases}
\end{gather}
We can eliminate the odd variables and get an equation for the even
variables. Let $u_n=x_{2n}$. Using $x_{2n-1}=u_nu_{n-1}-1$, the even
variables satisfy the recursion relation
\begin{equation*}
u_{n+1}={u_n^3+u_{n-1}\over u_n u_{n-1}-1},
\end{equation*}
or, equivalently,
$u_n(u_{n+1}u_{n-1}-u_n^2)=u_{n+1}+u_{n-1}$.

\subsubsection{Conserved quantities}
The variable $w_n=u_{n+1}u_{n-1}-u_n^2$ satisf\/ies
\begin{equation}\label{wn}
w_n={u_{n+1}+u_{n-1}\over u_n}={u_n^3+u_{n-1}+(u_nu_{n-1}-1)u_{n-1}\over u_n(u_nu_{ -1}-1)}
={u_n^2+u_{n-1}^2\over u_nu_{n-1}-1}.
\end{equation}
Moreover,
\begin{equation*}
w_{n+1}-w_n= {u_{n+1}^2+u_n^2\over u_{n+1}u_n-1}-{u_{n+1}+u_{n-1}\over u_n}=
{u_n^3-u_{n+1}u_nu_{n-1}+u_{n+1}+u_{n-1}\over u_n(u_nu_{n+1}-1)}=0.
\end{equation*}
We conclude that $w_n$ is a conserved quantity, that is, it is
independent of $n$.

Using \eqref{wn}, there exists a constant $c$ such that
$u_{n+1}+u_{n-1}=c u_n$. We may compute this constant explicitly in
terms of the initial conditions in cases $0$ and $1$ of \eqref{cases}, using $u_0=x_0$,
$u_1=x_2$ and $x_0x_2=1+x_1$:
\begin{gather*}
{\rm case}\ 0: \ \   c^{(0)} = {u_1^2+u_{0}^2\over u_1u_{0}-1}={x_0^2+x_2^2\over x_1}={x_0^4+(1+x_1)^2\over x_0^2x_1},\\
{\rm case}\ 1: \ \   c^{(1)} = {x_2^4 +(1+x_1)^2\over x_2^2 x_1}.
\end{gather*}

\subsubsection{Generating function}
The linear recursion relation $u_{n+1}-cu_n+u_{n-1}=0$
implies the following formulas for the generating functions
\begin{gather}
U^{(0)}(t) := \sum_{n\geq 0} u_n t^n={u_0-t(c^{(0)}u_0-u_1)\over 1-c^{(0)}t +t^2}
={x_0\over 1-t y^{(0)}_1 -{t^2 y^{(0)}_2\over 1-t y^{(0)}_3}},\label{casei} \\
U^{(1)}(t) := \sum_{n\geq 0}
u_{n+1} t^n={u_1-t(c^{(1)}u_1-u_2)\over 1-c^{(1)}t +t^2}
={x_2\over 1-t y_1^{(1)} -{t^2 y^{(1)}_2\over 1-t y^{(1)}_3}},\label{caseii}
\end{gather}
where the parameters are expressed in terms of the initial data as:
\begin{gather}
 y^{(0)}_1={1+x_1\over x_0^2},\qquad y^{(0)}_2={x_0^4+(1+x_1)^2\over x_0^4x_1},
\qquad y^{(0)}_3={x_0^4+1+x_1\over x_0^2x_1},\nonumber \\
 y^{(1)}_1={x_2^4+1+x_1\over x_2^2x_1},\qquad y^{(1)}_2={x_2^4+(1+x_1)^2\over x_2^4x_1},
\qquad y^{(1)}_3={1+x_1\over x_2^2}.\label{weightsonefour}
\end{gather}
Note that both sets satisfy the same relation
$y^{(i)}_1y^{(i)}_3=1+y^{(i)}_2$ where $i=0,1$,
as they are a~related via the substitution
$x_2\leftrightarrow x_0$, which maps $\big(y^{(1)}_1,y^{(1)}_2,y^{(1)}_3\big)\leftrightarrow
\big(y^{(0)}_3,y^{(0)}_2,y^{(0)}_1\big)$.

Upon expanding the right hand sides of
\eqref{casei} and \eqref{caseii} in power series of $t$,
we get all the $u_n$'s as explicit positive Laurent polynomials of
either initial data $(x_0,x_1)$ and $(x_1,x_2)$.  Using
\begin{equation*}
{1\over 1-t a_1-{t^2a_2 \over 1-t a_3}}=\sum_{p,q,\ell\geq 0} t^{p+2q+\ell} a_1^p a_2^q a_3^\ell
{p+q\choose p} {q+\ell-1\choose \ell}\end{equation*}
we obtain the expressions, valid for all $n\geq 0$:
\begin{gather}
\hbox{case $0$}: \ \     x_{2n}(x_0,x_1) =
\sum_{q,\ell,r,s,m\geq 0}
x_0^{1+4(q+\ell)-2n-4(r+s)} x_1^{m-q-\ell} \nonumber \\
\hphantom{\hbox{case $0$}: \ \     x_{2n}(x_0,x_1) = }{} \times
{n-q-\ell\choose q-r,r}{q+\ell-1\choose \ell-s,s}{n+2r+s-2q-\ell\choose m},
\label{solsonefourone} \\
\hbox{case $1$}: \ \ x_{2n+2}(x_1,x_2) =
\sum_{q,\ell,r,s,m\geq 0}
x_2^{1+2n-4(q+\ell+r+s)} x_1^{q+\ell+m-n} \nonumber \\
\phantom{\hbox{case $1$}: \ \ x_{2n+2}(x_1,x_2) =}{} \times {n-q-\ell\choose q-r,r,s}{q+\ell-1\choose \ell}{2r+s+\ell\choose m},
\label{solsonefourtwo}
\end{gather}
where we have used the multinomial coef\/f\/icients
\[
{n\choose m_1,m_2,\dots ,m_k}:=\left\{\begin{array}{ll} \dfrac{n!}{m_1! \cdots m_{k}!(n-\sum
  m_i)!}, &\sum_i m_i \leq n,\vspace{2mm}\\
0 & {\rm otherwise}.\end{array}\right.
\]

The expressions for the odd variables follow from
$x_{2n+1}=x_{2n}x_{2n+2}-1$. The positivity of these latter
expressions is easily checked: Note that the product $x_{2n}x_{2n+2}$
in both
equations~\eqref{solsonefourone} and~\eqref{solsonefourtwo}
has a constant term $1$ as a Laurent polynomial of
$(x_0,x_1)$ and $(x_1,x_2)$ respectively, hence
$x_{2n}x_{2n+2}-1$ is a positive Laurent polynomial.

\subsubsection{Path interpretation}
\begin{figure}[t]
\centerline{\includegraphics{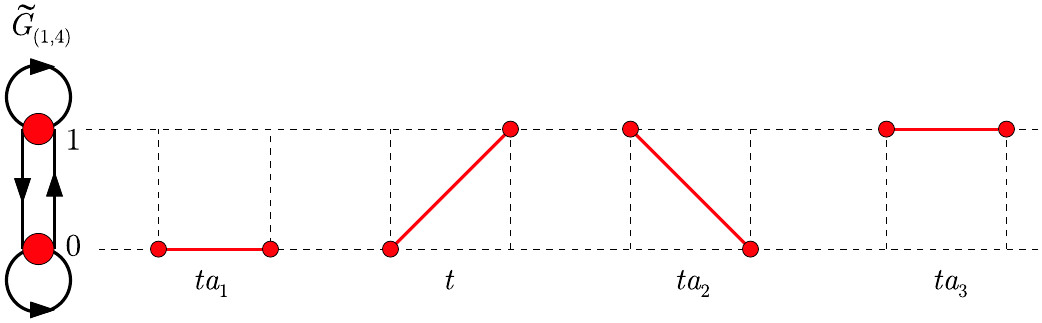}}
\caption{The graph ${\widetilde G}_{(1,4)}$ and the four corresponding path
steps.}\label{fig:graphonefour}
\end{figure}

The continued fraction expressions~\eqref{casei} and \eqref{caseii}
allow for a path interpretation of $x_{2n}$ as follows.  Consider the
graph ${\widetilde G}_{(1,4)}$ on the left hand side of Fig.~\ref{fig:graphonefour}, with two vertices labelled~$0$,~$1$ connected by
an edge, and connected to themselves via a loop.  We assign weights to
the oriented edges as follows:
\begin{gather*}
 w(0\to 0)=t a_1, \qquad  w(0\to 1)=t,\qquad
w(1\to1) =t a_3,\qquad w(1\to0) = t a_2.
\end{gather*}
We can also associate a path in $\Z_{\geq 0}^2$ to a path on
$\widetilde{G}_{(1,4)}$
composed of the steps shown in Fig.~\ref{fig:graphonefour}.

The corresponding path transfer matrix is:
$T=t \left(\begin{matrix}
    a_1  & a_2 \\
    1 & a_3
\end{matrix}\right)$.
Using Gaussian elimination on $I-T$ as in the previous example, we
f\/ind that the partition function of paths from the vertex~0 to itself
is:
\begin{equation*}
\big((I-T)^{-1}\big)_{0,0}= {1\over 1- t a_1 -{t^2 a_2\over 1-t a_3}}.
\end{equation*}

Hence, $x_{2n}$ is (up to a factor of $x_0$ in case $0$ and $x_2$ in case
$1$) the partition function for paths of $2n$ steps on
${\widetilde G}_{(1,4)}$, from and to the vertex $0$, with weights
$(a_1,a_2,a_3)=(y_1,y_2,y_3)$ in case $(i)$ and $(a_1,a_2,a_3)=(z_1,z_2,z_3)$ in
case $(ii)$, where the weights are as in \eqref{weightsonefour}.

\subsubsection[The graph $\widetilde{G}_{1,4}$ as a compactif\/ication of
  the graph $\widetilde{G}_1$]{The graph $\boldsymbol{\widetilde{G}_{1,4}}$ as a compactif\/ication of
  the graph $\boldsymbol{\widetilde{G}_1}$}\label{twopathone}

   We will see later
that, quite generally, it is possible to ``compactify'' the graphs
$\widetilde G_r$. The result for $r=1$ is precisely the graph
$\widetilde{G}_{1,4}$. Thus,
The case $(b,c)=(2,2)$ of the previous section may
also be interpreted in terms of paths on the graph
${\widetilde G}_{(1,4)}$, but with dif\/ferent weights:
\begin{gather*}
  w(0\to 0)=t y_1, \qquad  w(0\to 1)=1,\qquad
w(1\to1) =t y_3,\qquad w(1\to0) = t y_2.
\end{gather*}
where $y_i$ are as in equation \eqref{yweights}.

To see this, note that the continued fraction~\eqref{aoneres} may be
rewritten as:
\begin{gather*}
X(t)={x_0\over 1-t {y_1\over 1-t {y_2\over 1-t y_3}}}=x_0+ {tx_1\over 1-t y_1-{ty_2\over 1-t y_3}}=x_0+t x_1
\left(\left(I -\left( \begin{matrix}
t y_1 & ty_2\\
1 & t y_3
\end{matrix}\right)\right)^{-1}\right)_{0,0}.
\end{gather*}

Thus, for $n\geq 0$, the solution $x_{n+1}$ of the $(2,2)$ case is up
to a factor $x_1$ the partition function for paths on ${\widetilde
  G}_{(1,4)}$, from and to the origin vertex $0$, and with a total of
$n$ steps of the form $0\to 0$, $1\to 1$ or $1\to 0$, with respective
weights $y_1$, $y_2$, $y_3$.

\section[Application to the $A_r$ $Q$-system]{Application to the $\boldsymbol{A_r}$ $\boldsymbol{Q}$-system}\label{qsysappli}

In \cite{Ke07,DFK08}, we showed that the recursion relations
($Q$-systems) satisf\/ied by the characters of the Kirillov--Reshetikhin
modules of the quantum af\/f\/ine algebras associated with any simple Lie
algebra can be described in terms of mutations of a cluster
algebra. Solutions to the cluster algebra recursion relations are more
general, in that the initial conditions are not specialized, as they
are in the original $Q$-system satisf\/ied by characters of
Kirillov--Reshetikhin modules~\cite{KR}.

A particularly simple example of this is the case when the Lie algebra is
of type $A$. Characters of $A_n$ KR-modules, which are just Schur
functions corresponding to rectangular Young diagrams, are given by
the cluster variables upon specialization of the boundary conditions.
In~\cite{DFK08a}, we gave the general solution for the cluster
variables without specialization of initial conditions. For these
variables, we proved the positivity conjecture of Fomin and
Zelevinsky by mapping the problem to a partition function of weighted
paths on a graph. Let us review the results obtained in~\cite{DFK08a}.

\subsection{Def\/inition}

The $A_r$ $Q$-system is the following system of
recursion relations for a sequence $R_{\al,n}$, $\al\in
I_r=\{1,2,\dots ,r\}$ and $n\in \Z$:
\begin{equation}\label{qsys}
R_{\al,n+1}R_{\al,n-1}=R_{\al,n}^2+R_{\al+1,n}R_{\al-1,n},\end{equation}
with boundary conditions
\[ R_{0,n}=R_{r+1,n}=1\qquad
\hbox{for all} \ \ n\in\Z.
\]
We note that the case $r=1$ of $A_1$ coincides with the case $(2,2)$
treated in Section~\ref{twotwo}.

\begin{remark} The original $Q$-system, which is the one
satisf\/ied by the characters of KR-modules, dif\/fers from the system
\eqref{qsys} not only in that it is specialized to the initial
conditions $R_{\al,0}=\pm1$ for all $\al\in I_r$, but also by a minus
sign in the second term on the right. In this discussion, we choose to
renormalize the variables for simplicity. It is also possible (see the
appendix of \cite{DFK08}) to consider the recursion relation with
nontrivial coef\/f\/icients. In that case both \eqref{qsys} and the
original $Q$-system result from a specialization of the coef\/f\/icients.
\end{remark}

We wish to study the solutions of equation \eqref{qsys} for any given
initial data. Our standard initial data for the $Q$-system are the
variables $\bx_0=(R_{\al,0},R_{\al,1})_{\al\in I_r}$.  We may then
view equation~\eqref{qsys} as a three-term recursion relation in $n$,
which requires two successive values of $n$ (for each~$\al$) as
initial data. Our f\/irst goal is to express all $\{R_{\al,n}\}_{\al,n}$
in terms of the initial data~$\bx_0$.

\subsubsection{Symmetries of the system}

We may use equation \eqref{qsys} to get both $n\geq 2$ and $n\leq -1$
values in terms of $\bx_0$. These are related via the manifest
symmetry of the $Q$-system under $n\leftrightarrow 1-n$, which implies
$R_{\al,1-n}\big((R_{\al,0},R_{\al,1})_{\al\in I_r}\big)
=R_{\al,n}\big((R_{\al,1},R_{\al,0})_{\al\in I_r}\big)$.
In addition, equation~\eqref{qsys} is translationally invariant under
$n\to n+k$: $R_{\al,n+k}\big((R_{\al,k+1},R_{\al,k})_{\al\in I_r}\big)
=R_{\al,n}\big((R_{\al,1},R_{\al,0})_{\al\in I_r}\big)$ for all
$n,k\in\Z$.

\subsection[Conserved quantities of $Q$-systems as hard particle partition functions]{Conserved quantities of $\boldsymbol{Q}$-systems as hard particle partition functions}

The $Q$-system turns out to be a discrete integrable system, in that
it is possible to f\/ind a~suf\/f\/icient number of integrals of the motion,
as in the $r=1$ case.

First, equation \eqref{qsys} may be used to express
$R_{\al,n}$ $(\al\geq 2)$ as polynomials in $\{R_{1,n}\}$:
\begin{equation}\label{deteralpha}
R_{\al,n}=\det_{1\leq i,j \leq \al} \, (R_{1,n+i+j-\al-1}).
\end{equation}
This is proved by use of the Desnanot--Jacobi identity for the minors
of the $(\al+1)\times(\al+1)$ matrix $M$ with entries
$M_{i,j}=R_{1,n+i+j-\al}$, $i,j\in I_{\al+1}$.

The boundary condition $R_{r+1,n}=1$ (for all $n$) together with
 \eqref{deteralpha}
implies the equation of motion which determines $R_{1,n}$:
\[
\det_{1\leq i,j\leq r+1}
(R_{1,n+i+j-r-2})=1.
\]
From equation \eqref{qsys}, it is clear that $R_{r+2,n}=0$ for all $n$.
Hence there exists a linear
recursion relation of the form
\begin{equation}\label{recugen}
\sum_{m=0}^{r+1} (-1)^m c_{r+1-m} \, R_{1,n+m} =0,
\end{equation}
with $c_0=c_{r+1}=1$. The fact that $c_{r+1-m}$ do not depend on $n$
follows from the fact that each row in the matrix in the determinant
for $R_{r+2,n}$ is just a shift in $n$ of any other row.

The constants $c_p$, $p=1,2,\dots ,r$ are the $r$ integrals of motion of
the $A_r$ $Q$-system. They can be expressed explicitly in terms of the
$R_{1,n}$s as follows:
\begin{equation*}
c_p=\det_{1\leq i\leq r+1\atop 1\leq j\leq r+2,\ j\neq r+2-p}\, R_{1,n+i+j-2}
\end{equation*}
independently of $n$. These quantities are similar to those found in~\cite{RESH}
for the so-called Coxeter--Toda integrable systems.

By using simple determinant identities, we show in Theorem~3.5 of
\cite{DFK08a} that $c_p$ satisfy recursion relations which allow us to
identify them as the partition functions $Z_p^{G_r}$ of $p$ hard
particles on the graph $G_r$ of Fig.~\ref{fig:graphgr}, where the
weights given by:
\begin{gather}\label{weights}
y_{2\al-1,k}= {R_{\al-1,k}R_{\al,k+1}\over R_{\al,k}R_{\al-1,k+1}},\quad
1\leq \al\leq r+1, \qquad y_{2\al,k}={R_{\al-1,k}R_{\al+1,k+1}\over
R_{\al,k}R_{\al,k+1}}, \quad 1\leq \al\leq r.\!\!
\end{gather}
This is true for any $k$: The functions $Z_p^{G_r}$ are independent of
the choice of $k$. Thus, unless otherwise stated, $y_\al$ will stand
for $y_{\al,0}$ below.

\subsection[$Q$-system solutions and paths]{$\boldsymbol{Q}$-system solutions and paths}

The linear recursion relation \eqref{recugen} allows to compute the generating function
$R^{(r)}(t)=\sum\limits_{n=0}^\infty t^n R_{1,n}$ explicitly. Indeed,
$R^{(r)}(t) Z^{G_r}(-t\by)$
is a polynomial of degree $r$, and it is easy to see that
\[
R^{(r)}(t) Z^{G_r}(-t\by)=R_{1,0}Z^{G_r}(0,-ty_2,\dots ,-ty_{2r+1}).
\]

Using~\eqref{heaps}, we may interpret $R^{(r)}(t)/R_{1,0}$ as
the generating function ${\mathcal Z}_{0,0}^{{\widetilde G}_r}(t\by)$
for paths on~${\widetilde G}_r$ with the weights~\eqref{weights},
say for $k=0$. In other words, $R_{1,n}/R_{1,0}$ is the partition function
for weighted paths of $2n$ steps on ${\widetilde G}_r$ starting and ending
at the origin (or starting at $(0,0)$ and ending at $(2n,0)$ in the
two-dimensional representation).

Comparing the determinant formula for $R_{\al,n}$ \eqref{deteralpha}
and the LGV formula \eqref{lgv} for families of $\al$ non-intersecting
paths, we have \cite{DFK08a}
\begin{lemma} The quantity $R_{\al,n}/(R_{1,0})^\al$ is the partition
function ${\mathcal Z}_{\bs,\be}^{{\widetilde G}_r}$ for $\al$
non-intersecting weighted paths on ${\widetilde G}_r$ with starting
and ending points $s_i=(2i-2,0)$, $e_i=(2n+2\al-2i,0)$,
$i=1,2,\dots ,\al$.
\end{lemma}
\begin{proof}
It is clear that ${\mathcal Z}_{s_i,e_j}^{{\widetilde G}_r}$ is the
partition function for paths from $(2i-2,0)$ to $(2n+2\al-2j,0)$, which is
equal to that from $(0,0)$ to $(2(n+\al+1-i-j),0)$ by translational
invariance.  This allows to identify the two determinants
\eqref{deteralpha} and \eqref{lgv}, up to an overall factor of
$(R_{1,0})^\al$, and the result follows.
\end{proof}

\begin{figure}[t]
\centerline{\includegraphics[width=12.5cm]{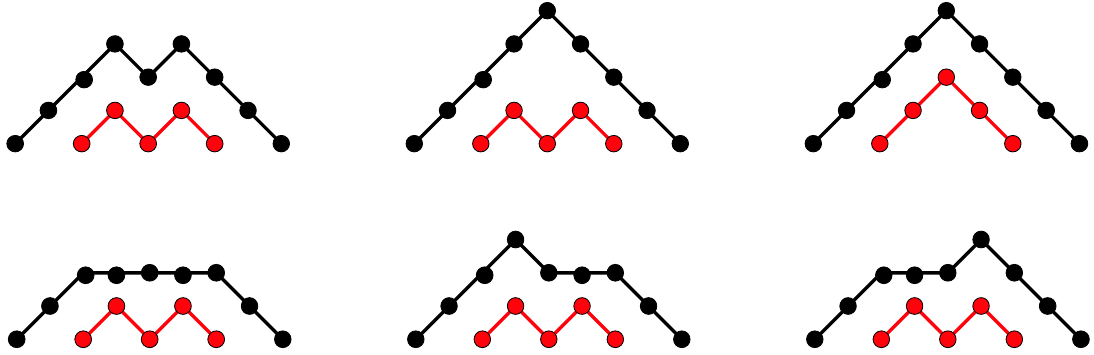}}
\caption{The six pairs of non-intersecting paths on ${\widetilde G}_2$ of $8$ and $4$ steps,
starting respectively at $(0,0)$ and $(2,0)$ and ending at
$(6,0)$ and $(8,0)$.}\label{fig:exatwo}
\end{figure}

 As an illustration, we have represented in Fig.~\ref{fig:exatwo} the
six pairs of non-intersecting paths contributing to $R_{2,3}$,
solution of the $A_2$ $Q$-system.

Thus, we have proved
\begin{theorem}
The variables $R_{\al,n}$ which satisfy \eqref{qsys}, when expressed
in terms of the va\-riab\-les~$\bx_0$, are equal to $(R_{1,0})^\al$ times
partition functions for paths on the graph~$\widetilde{G}_r$,
involving only the weights $y_\al$ of \eqref{weights}. These weights are
explicit Laurent monomials of the initial data $\bx_0=
(R_{\al,0},R_{\al,1})_{\al\in I_r}$. This gives an explicit expression
for $R_{\al,n}$ as Laurent polynomials of the initial data, with
non-negative integer coefficients.
\end{theorem}

\subsection{An alternative path formulation}\label{altersec}

\begin{figure}[t]
\centerline{\includegraphics{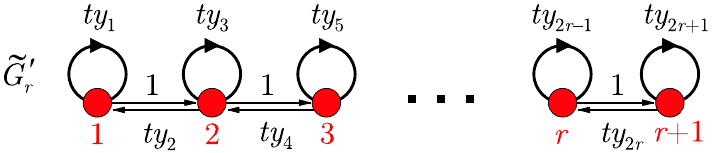}}
\caption{\small The graph ${\widetilde G}_r'$, with $r+1$ vertices. We have
indicated the weights attached to each oriented edge.}\label{fig:dualgrprime}
\end{figure}

In the same spirit as Section~\ref{twopathone}, one can show (see
below) that the solution $R_{1,n}$ of the $A_r$ $Q$-system may also be
interpreted in terms of paths on a new (``compactif\/ied'') graph
${\widetilde G}'_r$ of Fig.~\ref{fig:dualgrprime}. This is a graph
with $r+1$ vertices, labelled $1,2,\dots ,r+1$, connected via oriented
edges $i\to i+1$, $i+1\to i$, $i=1,2,\dots ,r$ and loops $i \to i$,
$i=1,2,\dots ,r+1$.

To each edge is attached a weight as follows:
\begin{gather}\label{newwts}
  {\rm wt}(i\to i+1) = 1, \qquad
{\rm wt}(i+1\to i) = ty_{2i},\qquad
{\rm wt}(i\to i)=t y_{2i-1}.
\end{gather}
The corresponding transfer matrix encoding these weights is an
$(r+1)\times (r+1)$-matrix of the form
\begin{equation}\label{transmatcompact}
T'= \left(
\begin{matrix}
ty_1 & ty_2 & 0        &  \cdots &  \cdots &  \cdots & \cdots & 0 \\
1 & ty_3     & ty_4 & 0    &   &   & &  \vdots\\
0 & 1     & ty_5    &  ty_6&  0 &    &  & \vdots \\
\vdots &  \ddots     &  \ddots    &   \ddots    &   \ddots & \ddots  &  & \vdots \\
\vdots &     &  \ddots     &   \ddots     & \ddots & \ddots & \ddots & \vdots \\
\vdots  &   &     &  0    &   1 & ty_{2r-3}            &  ty_{2r-2}         & 0 \\
\vdots  &  &     &       &   0 & 1         & ty_{2r-1}         & ty_{2r} \\
0 &   \cdots   &   \cdots   &    \cdots    &   \cdots & 0           & 1         & ty_{2r+1}
\end{matrix}
\right).
\end{equation}
Then we have
\begin{gather}
{R^{(r)}(t)\over R_{1,0}}=1+t{y_1\over 1-t y_1-t{y_2\over 1-t y_3 -t{y_4\over
1-ty_5-t{y_6\over {\ddots \over 1-ty_{2r-1}-t{y_{2r}\over 1-ty_{2r+1}}}}}}}\nonumber\\
\hphantom{{R^{(r)}(t)\over R_{1,0}}}{} =1+t {R_{1,1}\over R_{1,0}}
\big((I -T')^{-1}\big)_{1,1}. \label{newcont}
\end{gather}
This is readily proved by Gaussian elimination.

Comparing this expression to the results of the previous section,
we conclude that
\begin{lemma} For $n\geq 0$, $R_{1,n+1}/R_{1,1}$ is the partition
function for paths on the weighted graph~${\widetilde G}'_r$, with a
total of $n$ steps along the edges of type $i\to i$ or $i+1\to i$ in
${\widetilde G}'_r$, starting and ending at vertex $1$.
\end{lemma}
The weights
are given in equation \eqref{newwts} with $y_i$ as in \eqref{weights}
with $k=0$.

In Section~\ref{totalpos} we will relate this result to the total
positivity conjecture of Fomin and Zelevinsky and networks.

We may actually write an explicit expression for $R_{1,n}$ by simply expanding
the continued fraction \eqref{newcont} as:
\begin{gather*}
R^{(r)}(t)=R_{1,0} \Bigg(1+t y_1 \sum_{p_1,p_2,\dots ,p_{2r+1}\geq 0\atop p_0=p_{2r+2}=0}
\prod_{\ell=0}^{r} (ty_{2\ell+1})^{p_{2\ell+1}}\nonumber\\
\phantom{R^{(r)}(t)=}{} \times (ty_{2\ell+2})^{p_{2\ell+2}}
{p_{2\ell}+p_{2\ell+1}+p_{2\ell+2}-1\choose p_{2\ell}-1,p_{2\ell+1}}
\Bigg).
\end{gather*}
Substituting the values \eqref{weights} for the weights $y_{\al}\equiv y_{\al,0}$,
and extracting the coef\/f\/icient of $t^{n+1}$, we get for all $n\geq 0$:
\begin{gather*}
R_{1,n+1} = R_{1,1} \sum_{p_1,p_2,\dots ,p_{2r+1}\geq 0
\atop p_0=p_{2r+2}=0,\ \Sigma p_i=n}
\prod_{i=1}^{r} {(R_{i,0})^{p_{2i+2}+p_{2i+1}-p_{2i}-p_{2i-1}}\over
(R_{i,1})^{p_{2i+1}+p_{2i}-p_{2i-1}-p_{2i-2}}}
\prod_{\ell=0}^{r}{p_{2\ell}+p_{2\ell+1}+p_{2\ell+2}-1\choose p_{2\ell}-1,p_{2\ell+1}}
\end{gather*}
as explicit positive Laurent polynomials of the initial data. This gives a rank-$r$ generalization
of the formula given in \cite{CZ} for $r=1$.

\section[Cluster algebra formulation: mutations and paths for the $A_r$ $Q$-system]{Cluster algebra formulation: mutations and paths\\ for the $\boldsymbol{A_r}$ $\boldsymbol{Q}$-system}\label{generalpaths}

In this section, we show that the solutions $\{R_{\al,n}\ |\ \al\in
I_r,n\in \Z\}$ of the $Q$-system are positive Laurent polynomials when
expressed as functions of
an arbitrary set initial conditions. This generalizes our result
for the initial condition $\bx_0$ in the previous section.

The
recursion relation \eqref{qsys} has a solution once a certain set of
initial conditions is specif\/ied, but this set need not necessarily be
the set $\bx_0$. We will explain below that the most general possible
choice of initial conditions is specif\/ied by a Motzkin path of length $r$.

The solutions of \eqref{qsys} can be viewed as cluster variables in
the $A_r$ $Q$-system cluster algebra def\/ined in \cite{Ke07}. Hence,
our proof provides a general conf\/irmation of the conjecture of
\cite{FZI} for this particular cluster algebra: When the cluster
variables are expressed as functions of the variables in any other
cluster (i.e.\ an arbitrary set of initial conditions), they are
Laurent polynomials with non-negative coef\/f\/icients.

The results of this sections were explained in detailed in
\cite{DFK08a}, and this section should serve as a summary of the
proofs contained therein.

\subsection[The $Q$-system as cluster algebra]{The $\boldsymbol{Q}$-system as cluster algebra}

In \cite{Ke07}, it was shown that the $A_r$ $Q$-system solutions
$\{R_{\al,n}\}$ may be viewed as a subset of the cluster variables of
the $A_r$ $Q$-system cluster algebra.  This is a cluster algebra with
trivial coef\/f\/icients, which includes the the seed cluster variable
$(R_{1,0},R_{2,0},\dots ,R_{r,0},R_{1,1},R_{2,1},\dots ,R_{r,1})$, with an
associated associated $2r\times 2r$ exchange matrix has the block
form: $\begin{pmatrix}0 & -C \\ C & 0\end{pmatrix}$, $C$ the Cartan
matrix of $A_r$.

In this language, each cluster
is a vector with $2r$ variables, and the subset of clusters relevant
to the $Q$-system are those which have entries made up entirely of
solutions to the $Q$-system (we restrict to these in the following).
These clusters are all related by sequences of cluster
mutations which are one of the relations~\eqref{qsys}.

Because of the form of equation~\eqref{qsys}, it is easy to see that
the restricted set of clusters corresponding to the $Q$-system
are characterized by a set of $r$ integers
$(m_1,\dots ,m_r)$, subject to the condition that $|m_\al-m_{\al+1}|\leq
1$. This def\/ines what is known as a Motzkin path.
\begin{definition}\label{xmdef}
The cluster $\bx_{\bm}$ corresponding to the set of integers
$\bm=(m_1,\dots ,m_r)$ is the vector of $2r$ variables
$\{R_{\al,m_\al},R_{\al,m_\al+1}\}_{\al\in I_r}$, ordered so that all
variables with an even second index appear f\/irst.
\end{definition}
For example, the initial cluster $\bx_0$ corresponds to the Motzkin path
$\bm=(0,0,\dots ,0)$, and
$\bx_0=(R_{1,0},\dots ,R_{r,0},R_{1,1},\dots ,R_{r,1})$.

For any $\bm$, $\bx_\bm$ is obtained from
the fundamental initial {\em seed} $\bx_0$ by
{\em mutations} of the cluster algebra. This is just saying that one
gets $\bx_\bm$ by repeated selected applications of the recursion relation
\eqref{qsys} to $\bx_0$. Each mutation changes only one of the cluster
variables. That is, for some~$\al$ and~$n$,
\[
R_{\al,n}\mapsto\left\{
\begin{array}{ll}
R_{\al,n+2}=\frac{R_{\al,n+1}^2+R_{\al+1,n+1}R_{\al-1,n+1}}{R_{\al,n}}&
\hbox{(forward mutation at $\al$)},\\
R_{\al,n-2}=\frac{R_{\al,n-1}^2+R_{\al+1,n-1}R_{\al-1,n-1}}{R_{\al,n}}&\hbox
{(backward mutation at $\al$)}.\end{array}\right.
\]
Recall that we only consider here particular mutations that only involve solutions
of the $Q$-system.
Here we have used the ``time'' variable $n$ to def\/ine forward (resp. backward)
mutations according to whether the mutation increases (resp. decreases) the index $n$
in the mutated cluster variable.

Alternatively, the mutation changes the Motzkin path which
characterizes $\bx_\bm$, by changing  $\bm \mapsto \bm\pm
\epsilon_{\al}$, with the plus (minus) sign
for a forward (backward) mutation. Here, $\epsilon_{\al}$ is the
vector which is zero except for the entry $\al$, which is equal to~1.
We see that accordingly the Motzkin path locally moves forward (backward).

The Laurent property of cluster algebras ensures that
every cluster variable is a Laurent polynomial of the cluster
variables of any other cluster in the algebra.

The positivity property, proved only in particular cases so
far, is that these Laurent polynomials have non-negative integer coef\/f\/icients.
The property was proved for the particular clusters considered in the present case in~\cite{DFK08a}.
It may be stated as follows:
\begin{theorem}[\cite{DFK08a}]
Each
$R_{\al,n}$, when expressed as a function of the seed $\bx_\bm$
for any Motzkin path $\bm$, is a Laurent
polynomial of $\{R_{\al,m_\al},R_{\al,m_\al+1}\; | \; \al\in I_r\}$, with
non-negative integer coefficients.
\end{theorem}
We outline the proof below.

For clarity let us introduce the following notation. Let $F$ be some
cluster variable. Then~$F$ can be expressed as a function of $\bx_\bm$
for any $\bm$. The functional form is then denoted by
$F=F_\bm(\bx_\bm)$. Since $F$ can also be expressed as a function of
any other cluster, we can write $F_\bm(\bx_\bm) =
F_{\bm'}(\bx_{\bm'})$ for any two Motzkin paths $\bm$, $\bm'$. In
particular, in the notation of the previous section, we have
$R^{(r)}(t) = R^{(r)}_{\bm_0}(t;\bx_0)$.

\subsection{Target graphs and weights}

Due to the ref\/lection and translation symmetries of the $Q$-system,
we can restrict our attention to seeds associated with Motzkin paths
in a fundamental domain $\mathcal M_r=\{\bm \ | \ {\rm min}_\al( m_\al)
= 0\}$. There are $3^{r-1}$ elements in $\mathcal M_r$.

To each Motzkin path $\bm\in {\mathcal M}_r$, we associate a pair
$(\Gamma_\bm,\{y_e(\bm)\})$, consisting of a rooted graph $\Gamma_\bm$
with oriented edges, and edge weights $y_e({\bm})$ along the edges $e$.

\begin{figure}
\centerline{\includegraphics{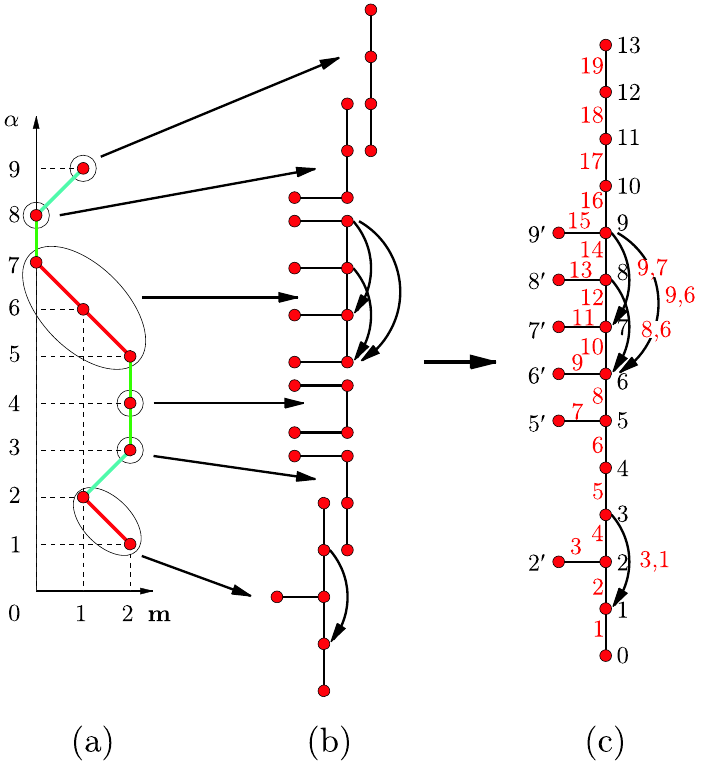}}
\caption{The Motzkin path
$\bm=(2,1,2,2,2,1,0,0,1)$ for $r=9$ (a) is decomposed into $p=6$ descending
  segments $(12)(3)(4)(567)(8)(9)$
(circled, red edges).
The corresponding graph pieces $\Gamma_{\bm_i}$ are indicated in~(b). They are to be glued ``horizontally'' for f\/lat
transitions (green edges) and ``vertically'' for ascending ones (blue edges).
The resulting graph $\Gamma_\bm$ is represented (c) with its
vertex (black) and edge (red) labels.}\label{fig:exglue}
\end{figure}

\subsubsection[Construction of the graph $\Gamma_\bm$]{Construction of the graph $\boldsymbol{\Gamma_\bm}$}

The graph $\Gamma_\bm$ is constructed via the following sequence of
steps (see Fig.~\ref{fig:exglue} for an illustration):
\begin{enumerate}\itemsep=0pt
\item Decompose the Motzkin path $\bm$ into maximal ``descending
segments'' $\bm_i$ of length $k_i$ ($i=1,\dots ,p$). These are segments of
the form $\bm_i=(m_{\al_i},m_{\al_i+1},\dots ,m_{\al_{i+1}-1})$ with
$\alpha_{i+1}=\alpha_i+k_i$, where $m_{\al_i+j}=m_{\al_i}-j$.
Here, $\al_1=1$ and $\al_{p+1}-1=r$.
\item The {\em separation} between two consecutive descending segments of
  the Motzkin path, $\bm_i$ and $\bm_{i+1}$ is either ``f\/lat''
  i.e.\ $m_{\al_{i+1}}=m_{\al_{i+1}-1}$ or ``ascending'' i.e.\
  $m_{\al_{i+1}}=m_{\al_{i+1}-1}+1$.
\item To each descending segment $\bm_i$, associate a graph
  $\Gamma_{\bm_i}$, which is the graph ${\widetilde G}_{k_i}$ with
  additional, down-pointing edges $a\to b$ for all $a$, $b$ such that
  $k_i+1\geq a>b+1\geq 2$. There are a total of $k_i(k_i-1)/2$ extra
  oriented edges.
\item We glue the graphs
$\Gamma_{\bm_i}$ and $\Gamma_{\bm_{i+1}}$ into a graph
$\Gamma_{\bm_i}\vert\vert \Gamma_{\bm_{i+1}}$ def\/ined as follows (see
Fig.~\ref{fig:exglue} for an illustration):
\begin{enumerate}\itemsep=0pt
\item If the separation between $\bm_i$ and $\bm_{i+1}$ is f\/lat, we
  identify vertex $0$ of $\Gamma_{\bm_{i+1}}$ with vertex $k_i+2$ of
  $\Gamma_{\bm_i}$, and vertex $1$ of $\Gamma_{\bm_{i+1}}$ with vertex
  $k_i+1$ of $\Gamma_{\bm_i}$, while the connecting edges are identif\/ied.
\item If the separation is ascending, we reverse the role of vertices
  $0$ and $1$ in the procedure above.
\end{enumerate}
\end{enumerate}
The result of this procedure is the graph
$\Gamma_\bm=\Gamma_{\bm_1}||\Gamma_{\bm_2}|| \cdots
||\Gamma_{\bm_p}$. Its root is the vertex~$0$ of~$\Gamma_{\bm_1}$.

We label the vertices of the graph $\Gamma_\bm$ by the integers $i$
with $i \in \{0,\dots,r+2+n_+(\bm)\}$ (where $n_+(\bm)$ is the number of $m_\al$
such that $m_{\al+1}=m_\al+1$) and labels $i'$ for any univalent vertex
attached to vertex $i$ via a horizontal edge. We do this
by labeling the vertices of
$\Gamma_\bm$ from bottom to top, by shifting the labels of the subgraphs
$\Gamma_{\bm_i}$ so that no label is skipped nor repeated.

The edges $e$ pointing towards the root of $\Gamma_\bm$ are of two types:
\begin{enumerate}\itemsep=0pt
\item[$(i)$] the
``skeleton edges'' belonging to some ${\widetilde G}_{k_i}$ in
the above construction;
\item[$(ii)$] the extra, down-pointing
edges added in the gluing procedure.
\end{enumerate}

\subsubsection[The weights on the graph $\Gamma_\bm$]{The weights on the graph $\boldsymbol{\Gamma_\bm}$}

We label the $2r+1$ skeleton edges of type~$(i)$
by $\al=1,2,\dots ,2r+1$ from bottom to top (see the example in
Fig.~\ref{fig:exglue}), and the weights are denoted by~$y_\al(\bm)$. Weights assigned to edges pointing away from the root
are all set to~1.

Alternatively, we may label the ``down pointing'' skeleton edges by
the pairs of vertices $i+1\to i$ or $i'\to i$ which they connect.  The
extra edges of type~$(ii)$ are also labeled by the pairs $a\to b$ of
vertices which they connect. All edge weights may be labeled by the label
of the edge.

The weights of the edges of type $(ii)$ can be
expressed in terms of the skeleton weights:
\begin{equation*}
y_{a,b}(\bm)={\prod\limits_{b\leq i <a} y_{i+1,i}(\bm) \over
\prod\limits_{b<i <a} y_{i',i}(\bm)},
\end{equation*}
so that they obey the following intertwining condition
\begin{equation}\label{intertw}
y_{a,b}(\bm)y_{a',b'}(\bm)=y_{a,b'}(\bm)y_{a',b}(\bm),\qquad a>a'>b>b'.
\end{equation}

For example,
the extra weights of the example of Fig.~\ref{fig:exglue}
read respectively: $y_{3,1}=y_2y_4/y_3$, $y_{9,7}=y_{14}y_{12}/y_{13}$,
$y_{8,6}=y_{12}y_{10}/y_{11}$, and
$y_{9,6}=y_{9,7}y_{8,6}/y_{12}=y_{14}y_{12}y_{10}/(y_{13}y_{11})$.

Finally, for a given Motzkin path $\bm\in{\mathcal M}_r$, we def\/ine
the skeleton weights $y_\al(\bm)$, $\al=1,2,\dots ,2r+1$ to be:
\begin{gather}
y_{2\al-1}(\bm) =  {\lambda_{\al,m_\al}\over \lambda_{\al-1,m_{\al-1}}} , \qquad \al=1,2,\dots ,r+1,
\label{oddy}\\
y_{2\al}(\bm) =  {\mu_{\al+1,m_\al+1}\over \mu_{\al,m_{\al}}}
\left(1-\delta_{m_\al,m_{\al+1}+1} +
{\lambda_{\al+1,m_{\al+1}}\over \lambda_{\al+1,m_{\al}}}\delta_{m_\al,m_{\al+1}+1}\right)
 \nonumber\\
\phantom{y_{2\al}(\bm) =}{}  \times
\left(1-\delta_{m_{\al-1},m_{\al}+1} +
{\lambda_{\al-1,m_{\al}}\over \lambda_{\al-1,m_{\al-1}}}\delta_{m_{\al-1},m_{\al}+1}\right) ,
\qquad \al=1,2,\dots ,r, \label{eveny}
\end{gather}
where
\begin{equation*}
\lambda_{\al,n}={R_{\al,n+1}\over R_{\al,n}} , \qquad \mu_{\al,n}={R_{\al,n}\over R_{\al-1,n}} .
\end{equation*}
Note that with these def\/initions the expressions
\eqref{oddy}, \eqref{eveny} involve only variables of the seed~$\bx_\bm$.

To each Motzkin path $\bm\in {\mathcal M}_r$, we may f\/inally associate
a transfer matrix $T_\bm\equiv T_\bm(t \by(\bm))$, with entries
$(T_\bm)_{b,a}=$ weight of the oriented edge $a\to b$ on
$\Gamma_\bm$. Then the series in $t$
\begin{equation}\label{zmgen}
\mathcal Z_\bm(t;\bx_\bm):=\big((I-T_{\bm})^{-1}\big)_{0,0}=\sum_n t^n \mathcal
Z_{0,0}^{\Gamma_{\bm}}(n)
\end{equation}
 is the generating function for weighted paths on $\Gamma_\bm$, with
the coef\/f\/icient of $t^n$ being the partition function of walks from
vertex 0 to itself on $\Gamma_\bm$ which have $n$ down-pointing
steps. When $\bm=\bm_0=\mathbf 0$, this coincides with~ \eqref{genpath}.

\subsection{Mutations, paths and continued fraction rearrangements}\label{mutasec}

Our purpose is to write an explicit expression for the functional
dependence of the variab\-les~$R_{\al,n}$ on the seed variable
$\bx_\bm$, that is, f\/ind $(R_{\al,n})_{\bm}(\bx_\bm)$ for each $\al$, $n$
and a Motzkin path $\bm$.

To do this, we will describe how the generating function
$R^{(r)}(t)=R^{(r)}_{\bm}(t;\bx_{\bm})$ is related to the generating
function $R^{(r)}_{\bm'}(t;\bx_{\bm'})$, where $\bm$ and $\bm'$ are
related by a mutation. Then we start from the known function
$R_{\bm_0}^{(r)}(t;\bx_{\bm_0})$, and apply mutations to obtain all
other functions $R^{(r)}_{\bm}$ with $\bm$ in the fundamental domain~$\mathcal M_r$.

One can cover the entire fundamental domain $\mathcal M_r$
starting from $\bm_0=\mathbf 0$ by using only forward mutations
$\bm\mapsto \bm'=\bm+\epsilon_\al$ of either type
$(i)$ $(\dots ,a,a,a+1,\dots )\mapsto(\dots ,a,a+1,a+1,\dots )$ and $(ii)$
$(\dots ,a,a,a,\dots )\mapsto (\dots ,a,a+1,a,\dots )$, with the obvious
truncations
when $\al=1$ or $r$. (See Remark~8.1 in \cite{DFK08a}.)

Suppose $\bm$ and $\bm'$ are related by such a mutation.
We compute the two generating
functions of the type \eqref{zmgen}, $\mathcal Z_\bm(t;\bx_\bm)$ and
$\mathcal Z_{\bm'}(t;\bx_{\bm'})$.

In fact, the two matrices, $T_{\bm}$ and $T_{\bm'}$ dif\/fer only
locally, so that in computing the two generating functions by row
reduction, we f\/ind that the calculation dif\/fers only in a f\/inite
number of steps. Note that generating functions take the form of
f\/inite continued fractions with manifestly positive series
expansions of $t$.

We note two simple rearrangement lemmas which can be used to relate
f\/inite continued fractions:{\samepage
\begin{gather*}
{(R_1)} \quad { 1\over  1-{ a\over 1-b}}=1+{a\over  1-a-b},\nonumber \\
{(R_2)} \quad  a+{b\over 1-c}={a'\over  1-{b'\over 1-c'}},\quad  \hbox{where} \ \
 a'=a+b,\  \ b'={bc\over a+b},\  \ c'={ac\over a+b}.
\end{gather*}
One checks this by explicit calculation.}

Let $\al>1$.
Then, using  $(R_2)$ we can show that $\mathcal
Z_\bm(t;\bx_\bm)= \mathcal Z_{\bm'}(t;\bx_{\bm'})$ if
and only if
the weights $y_i\equiv y_i(\bm)$ and $y_i'\equiv
y_i(\bm')$ are related via:
\begin{gather*}
(i) \ m_{\al-1}=m_\al < m_{\al+1}: \quad y_\beta'=\left\{\begin{array}{ll}
%
y_{2\al-1}+y_{2\al}, & \beta=2\al-1, \\
y_{2\al}y_{2\al+1}/(y_{2\al-1}+y_{2\al}), & \beta=2\al,\\
y_{2\al-1}y_{2\al+1}/(y_{2\al-1}+y_{2\al}), & \beta=2\al+1,\\
y_\beta, & \hbox{otherwise},\end{array}\right.\\
(ii)\ m_{\al-1}=m_\al = m_{\al+1}:  \quad y_\beta'=
\left\{\begin{array}{ll}
y_{2\al-1}+y_{2\al}, & \beta=2\al-1, \\
y_{2\al}y_{2\al+1}/(y_{2\al-1}+y_{2\al}), & \beta=2\al, \\
y_{2\al-1}y_{2\al+1}/(y_{2\al-1}+y_{2\al}), & \beta=2\al+1,\\
y_{2\al+2}y_{2\al-1}/(y_{2\al-1}+y_{2\al}), &\beta=2\al+2,\\
y_\beta, & \hbox{otherwise}.\end{array}\right.
\end{gather*}
One checks directly that the expressions
\eqref{oddy}, \eqref{eveny} indeed satisfy the above relations.

The boundary case, where $\al=1$, is treated analogously, but f\/irst
requires a ``rerooting" of the graph to its vertex $1$, which is
implemented by the application of $(R_1)$: Indeed, we simply write
$\mathcal
Z_{\bm}(t;\bx_\bm)=\big((I-T_\bm)^{-1}\big)_{0,0}=1+ty_1(\bm) \mathcal
Z'_\bm(t;\bx_{\bm})$ with $\mathcal
Z_\bm'(t;\bx_\bm)=\big((I-T_\bm)^{-1}\big)_{1,1}$. We then rearrange
$\mathcal Z_\bm'(t;\bx_{\bm})$ using $(R_2)$ again, and f\/ind that $\mathcal
Z_{\bm}'(t;\bx_\bm)=\mathcal Z_{\bm'}(t;\bx_{\bm'})$ if and only if
the weights are related via the above equations.

The net result is the following. Given a compound mutation $\mu_\bm$
which maps the fundamental Motzkin path $\bm_0$ to
$\bm=(m_\al)_{\al\in I_r}$, then there are exactly
$m_1$ ``rerootings''
as described above. This corresponds to rewriting
the generating function
\begin{gather*}
R^{(r)}_\bm(t;\bx_\bm)=\sum_{i=0}^{m_1-1}
t^iR_{1,i}+t^{m_1}R_{1,m_1}\cZ_{\bm}(t;\bx_\bm) 
\end{gather*}
with $\cZ_{\bm}(t;\bx_\bm) $ as in equation \eqref{zmgen}.

\begin{figure}[t]
\centerline{\includegraphics{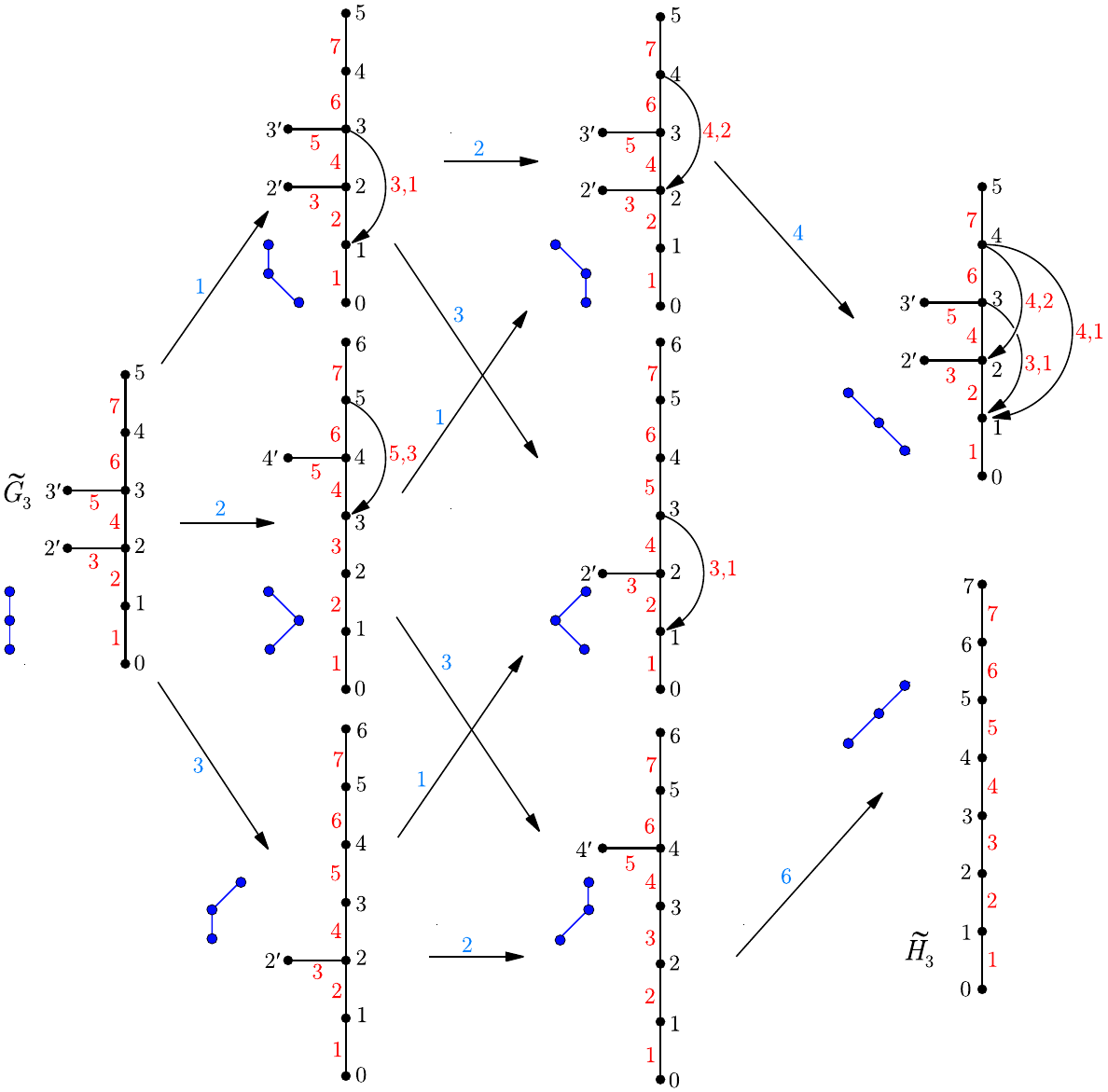}}
\caption{The Motzkin paths $\bm$ of the fundamental domain ${\mathcal M}_3$
and the associated graphs $\Gamma_\bm$, with their vertex and edge labels.
We have also indicated the
mutations by arrows, the label being $\al$ when the mutation $\mu_\al$ acts on variables $R_{\al,m}$
with an even index $m$
and $\al+r$ for an odd index $m$.}\label{fig:pathmotzfour}
\end{figure}

This leads to the following main result:

\begin{theorem}
For each $n\geq 0$, the function $(R_{1,n+m_1})_\bm(\bx_\bm)=
R_{1,m_1} \, {\mathcal Z}_{0,0}^{\Gamma_\bm}(n)$.  Thus it is
proportional to the generating function for weighted paths on the
graph $\Gamma_\bm$ with positive weights, so it is a manifestly
positive Laurent polynomial of the initial data $\bx_\bm$.
\end{theorem}

We have represented in Fig.~\ref{fig:pathmotzfour} the graphs $\Gamma_\bm$
for the Motzkin paths $\bm$ of the fundamental domain ${\mathcal M}_r$ for $r=3$.

\subsection[$Q$-system solutions as strongly non-intersecting paths]{$\boldsymbol{Q}$-system solutions as strongly non-intersecting paths}

To treat the case of $R_{\al,n}$ with $\al\neq 1$, given a Motzkin
path $\bm\in {\mathcal M}_r$, we need a path interpretation for the
determinant formula for $R_{\al,n+m_1}$:
\begin{gather} \label{deteralgv}
{R_{\al,n+m_1}\over (R_{1,m_1})^\al}= {\det \left(
R_{1,n+m_1+i+j-\al-1}\right)_{1\leq i,j \leq \al}\over
R_{1,n+m_1}}=\det\left( {\mathcal
Z}_{0,0}^{\Gamma_\bm}(n+i+j-\al-1)\right)_{1\leq i,j \leq \al}.
\end{gather}
Here, we have used the result of the previous section to rewrite the
formula in terms of the partition function for paths on $\Gamma_\bm$,
from and to the root, and with $n+i+j-\al-1$ down steps. As in the
standard LGV formula, we interpret this determinant as a certain
partition function for paths on $\Gamma_\bm$ starting from the root at
times $0,2,4,\dots ,2\al-2$ and ending at the origin at times
$2n,2n+2,\dots ,2n+2\al-2$.

\begin{figure}[t]
\centerline{\includegraphics{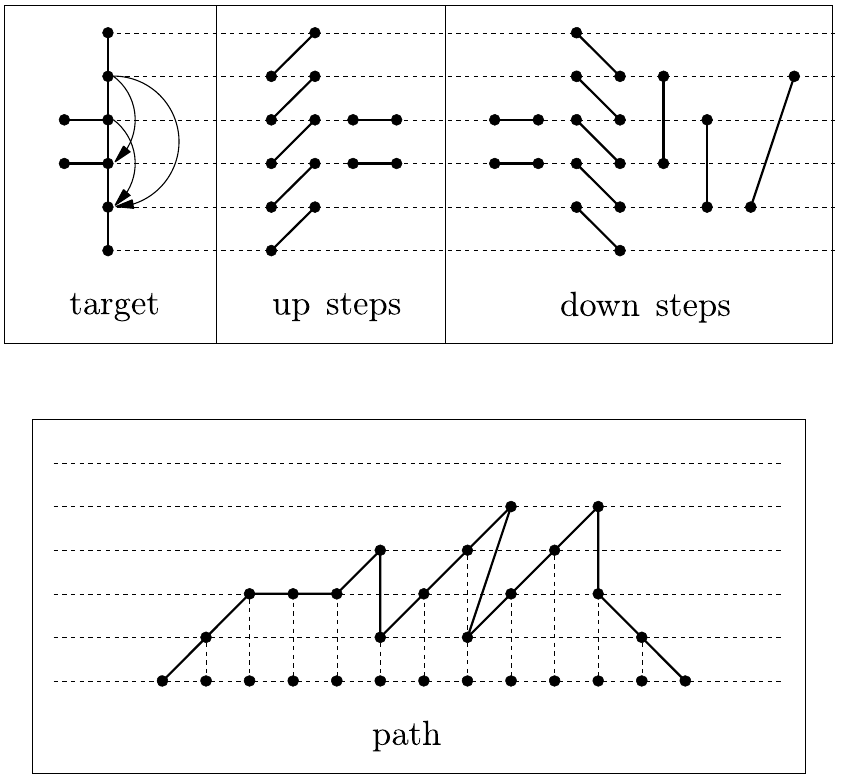}}
\caption{The two-dimensional representation of a typical path on the graph
$\Gamma_\bm$, $\bm$ the strictly descending Motzkin path $(2,1,0)$
of the case $A_3$.
Descents of $h=2$ are vertical (time displacement by $2-h=0$), while
descents of $h=3$ go back one step in time  (time displacement by $2-h=-1$). With these
choices, the total time distance between start and end is twice the number of descents
($16=2\times 8$ here).}\label{fig:pathslope}
\end{figure}

\subsubsection[Paths on $\Gamma_\bm$ represented as paths on a square lattice]{Paths on $\boldsymbol{\Gamma_\bm}$ represented as paths on a square lattice}

We draw paths, with allowed steps dictated by the graph $\Gamma_\bm$, on a
square-lattice in two dimensions. Paths start and end at
$y$-coordinate $0$. Moreover, if a path has $n$ ``down'' steps (steps
towards the vertex 0), then its starting and ending point are
separated by $2n$ horizontal steps. That is, a path is from $(x,0)$ to
$(x+2n,0)$ where $x$ the starting time and $n$ is the number of down-steps.

Since the horizontal distance between the starting and ending points
is f\/ixed by the {\em number} of ``down'' steps, a single step of the form
\[
a\to b=a-h
\]
should have a horizontal displacement $2-h$ (instead of 1 as in the
usual case). That is, on the square lattice it is a segment of the form
\[
(x,a) \to (x+2-h,a-h).
\]
Some examples are illustrated in Fig.~\ref{fig:pathslope}.

Thus, we identify ${\mathcal
Z}_{0,0}^{\Gamma_\bm}(n+i+j-\al-1)={\mathcal
Z}_{s_i,e_j}^{\Gamma_\bm}$ with paths on the two-dimensional lattice
starting at the point $s_i=(2i-2,0)$ and ending at the point
$e_j=(2n+2\al-2j,0)$, with the types of steps allowed given by the
edges of $\Gamma_\bm$ in the way explained in the previous paragraph.

\subsubsection{Strongly non-intersecting paths}
We now look at families of $\al$ paths, corresponding to the
determinant in equation \eqref{deteralgv}. Such paths may have
crossing on the lattice.  As in the case of LGV formula, the
determinant cancels out contributions from paths which share a
vertex. However, other situations may occur: Two paths may cross
without sharing a vertex in our picture.

One can generalize the proof of the LGV formula to take such crossings
into account. Using the expansion \eqref{expandet}, and introducing an
involution $\varphi$ on families of paths. This involution
interchanges the beginnings of the f\/irst two paths which share a
vertex or which cross each other, by transforming the crossing
segments $[P,Q]$ and $[R,S]$ into non-crossing ones $[R,Q]$ and
$[P,S]$. This ef\/fectively interchanges the two paths up to the points
$P$ and $R$ respectively, whichever comes f\/irst.

As in the usual case, the involution $\varphi$ acts as the identity if
no two paths cross, share a vertex, or can be made to cross via such
an exchange.

Taking into account the weights of the paths, the
intertwining condition \eqref{intertw} implies that the f\/lip preserves
the absolute value of the weight, but changes its sign, due to the
transposition of starting points. So the determinant \eqref{deteralgv}
cancels not only the paths that share a vertex or that cross, but also
those that come ``too close'' to one-another, namely that can be made
to cross via a f\/lip.

We call the families of paths which are invariant under the involution
$\varphi$ strongly non-intersecting paths.

To summarize, we have the following theorem:
\begin{theorem}
For any Motzkin path $\bm\in \mathcal M_r$, the variable
$R_{\al,n+m_1}$ $($with $n+m_1\geq \al-1)$ is equal to $(R_{1,m_1})^\al
$ times the partition function of $\al$ strongly non-intersecting
paths with steps and weights determined by $\Gamma_\bm$. The starting
points are $s_i=(2i-2,0)$ and the end points are $e_j=(2n+2\al-2j,0)$,
with $i,j=1,\dots ,\al$, and the weights are functions of the cluster~$\bx_\bm$.
\end{theorem}

In particular, $(R_{\al,n+m_1})_\bm(\bx_\bm)$ is a Laurent polynomial
with non-negative coef\/f\/icients of the cluster $\bx_\bm$.

\section[A new path formulation for the $A_r$ $Q$-system]{A new path formulation for the $\boldsymbol{A_r}$ $\boldsymbol{Q}$-system}\label{section6}

In Section \ref{generalpaths}, we constructed a set of transfer
matrices $T_\bm$, associated with paths on the graphs~$\Gamma_\bm$,
which allowed us to interpret~$R_{1,n}$ as generating functions of
weighted paths on a~graph~$\Gamma_\bm$, and hence prove their
positivity as a function of the seed variables $\bx_\bm$.

In Section \ref{altersec}, we also showed that for the special case
$\bm=\bm_0=\mathbf 0$, there is an alternative graph
$\Gamma'_{\bm_0}=\widetilde{G}_r'$, and that one can interpret
$R_{1,n}$ as a generating function for paths this alternative
``compactif\/ied'' graph. The graph $\Gamma'_{\bm_0}$ has $r+1$
vertices, hence the associated transfer matrix is of size $r+1\times
r+1$.

We now ask the question, is there a corresponding compactif\/ied set of
graphs, $\Gamma_\bm'$, which give a path formulation of $R_{1,n}$ with
weights which given by functions of the seed $\bx_\bm$ for all Motzkin
paths in~$\mathcal M_r$?

It turns out that it is always possible to f\/ind a weighted graph with
$r+1$ vertices which answers this question positively for each $\bm$. This
corresponds to a set of transfer matrices of a~size equal to
the rank of the algebra $A_r$. Therefore this transfer matrix
allows us to make a~direct connection between our transfer matrix
approach and the totally positive matrices of~\cite{FZposit}.

\subsection{Compactif\/ied graphs}
Consider the collection of graphs $\Gamma_\bm$. If we are interested
in the generating function of weighted paths on them from the vertex
$0$ to itself, then we can make various changes in them locally
(``compactify'' them) without af\/fecting the generating function
itself.

We obtain such compactif\/ied graphs from $\Gamma_\bm$ by identifying
pairs of neighboring vertices, and, when necessary, adding oriented
edges to cancel unwanted terms. There are two possible ways to make
such identif\/ications. Before presenting the general case, let us
illustrate the two situations in the following subsection.

\begin{figure}[t]
\centerline{\includegraphics{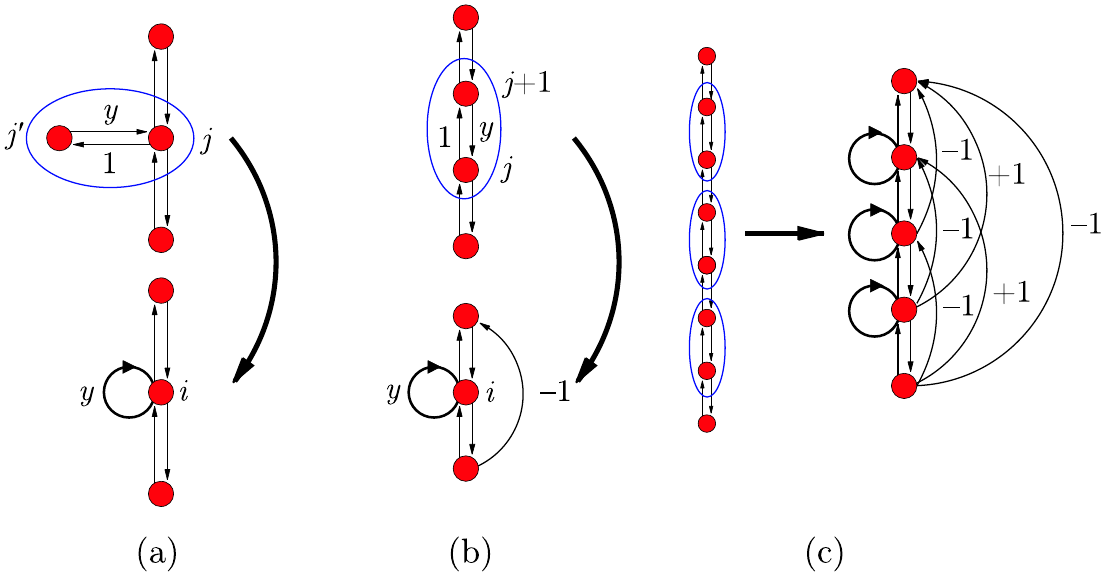}}
\caption{The identif\/ication of horizontal (a) or vertical (b) pairs of
consecutive vertices on $\Gamma_\bm$ and the result on $\Gamma_\bm'$.
The new edge on $\Gamma_\bm'$ with weight $-1$ allows to subtract the contribution from paths
that do not exist on $\Gamma_\bm$. We have also represented the situation
of a longer chain (c), where the identif\/ication now requires a network of up-pointing
edges with alternating weights $\pm 1$ for the suitable
subtractions.}\label{fig:loopsit}
\end{figure}

\subsection{Examples of compactif\/ication}
In fact, the f\/irst type of compactif\/ication, applied to
$\widetilde{G}_r$, leads to the alternative graph~$\widetilde{G}_r'$
obtained in Section~\ref{altersec}.  Recall that the generating
function for paths from $1$ to $1$ on the graph~$\widetilde{G}_r$ of
Fig.~\ref{fig:dualgrter} is related to the generating function for
paths from $0$ to $0$ via the rerooting procedure:
\[
\cZ_{0,0}^{\widetilde{G}_r}=1+ty_1\times \cZ_{1,1}^{\widetilde{G}_r}.
\]

\begin{example}
The generating function $\cZ_{1,1}^{\widetilde{G}_r}$
is equal to the generating function of paths from~$1$ to~$1$ on the graph
obtained from $\widetilde{G}_r$ via the
following procedure (``compactif\/ication''):
\begin{enumerate}\itemsep=0pt
\item Identifying vertex $i'$ with vertex $i$, whenever both exist, and
attaching a loop with weight $t y_{i',i}$ to the resulting vertex (see
Fig.~\ref{fig:loopsit}~(a)).
\item Identifying vertex $r+2$ with vertex $r+1$, and attaching a loop
  with weight $t y_{r+2,r+1}$ to the resulting vertex.
\item Identifying vertex $1$ with vertex $0$, renaming the resulting
  vertex $1$, and attaching a loop with weight $t y_1=t y_{1,0}$ to
  this vertex.
\end{enumerate}
This is clear, since a path from~$1$ to $1$ through the vertex $i'$
must have a segment $i\to i' \to i$ with weight $t y_{i',i}$, and if
it goes through the vertex $r+2$ it must have a segment $r+1\to r+2
\to r+1$ with weight $t y_{r+2,r+1}$. Similarly, the loop at $1$ accounts for
segments of the form $1\to 0 \to 1$.

Note that the resulting graph is $\widetilde{G}_r'$ of Fig.~\ref{fig:dualgrprime}. Thus, we have the identity of Section~\ref{altersec}
\[
\cZ_{1,1}^{\widetilde{G}_r} = \cZ_{1,1}^{\widetilde{G}_r'}.
\]
\end{example}

Another example of identif\/ication of vertices is the case when the two
vertices are adjacent vertices, $(i,i+1)$, on the spine of $\Gamma$
(see Fig.~\ref{fig:loopsit} (b)).

\begin{example}
Consider the graph $\widetilde{H}_k$ associated to an ascending
Motzkin path segment of length $k$ (see the example for $k=3$ in the
lower right hand corner of Fig.~\ref{fig:pathmotzfour}).  This is a
vertical chain of $2k+2$ vertices, numbered from $0$ to $2k+1$ from
bottom to top, connected by edges oriented in both directions. The
edges $i+1\to i$ have weights $y_{i+1}$ and the edges $i\to i+1$ have
weights~1.

Suppose we identify the vertices $a,a+1$ in this graph for some $a$. A
path from $0$ to $0$ with a step $a\to a+1$ is always paired with a
step $a+1\to a$, for a net contribution to the weight of the path of
$y_{a+1}$. We associate a loop with weight $y_{a+1}$ at the newly
formed vertex after the identif\/ication of vertices $a$ and $a+1$.

However, there are ``forbidden'' paths on the resulting graph, paths
which are not inherited from paths on $\widetilde{H}_r$. These are
paths which go from $a+2\to a-1$ without traversing the loop. This
would correspond to going from $a+2$ to $a-1$ on $\widetilde{H}_r$
without passing through the edge $a+1\to a$, which is impossible.
We cancel the contribution of these paths by adding an ascending edge
$a-1\to a+2$ with weight $-1$. The ef\/fect is precisely to subtract the
weights of the forbidden set of paths.
\end{example}

More generally, a succession of identif\/ications of the type in the
example above results in the following (see Fig.~\ref{fig:loopsit}~(c)):
\begin{lemma}\label{inclexcl}
The generating function of paths from $0$ to $0$ on the graph
$\widetilde{H}_k$ with vertices $0,\dots ,2k+1$ is equal to
the generating function for paths from $0$ to $0$ on the following
compactified graph~$\w H'_k$:
\begin{enumerate}\itemsep=0pt
\item[$1.$] Identify the vertices $2i+1$ and $2i+2$; rename the resulting
  vertex $i+1$ $(i=0,\dots ,k-1)$.
\item[$2.$] Attach a loop at the vertex $i+1$ with weight
  $w_{i+1}=y_{2i+2}$. Other edge weights remain unchanged.
\item[$3.$] Add ascending edges $j\to j+2+a$ $(0\leq a\leq k-1-j$, $0\leq j \leq k-1)$
  with weight $(-1)^{a+1}$ to the resulting graph.
\end{enumerate}
\end{lemma}
An illustration of the resulting graph $\w H'_k$ is given in (c) of
Fig.~\ref{fig:loopsit}.

\begin{proof}
Consider the set of paths $\mathcal P$ of the form $P_1 P^+ P_2 P^-
P_3$, where $P_i$ are f\/ixed paths, and~$P^+$ is a path from~$h_0$ to
$h_1$ consisting of only up steps and loop steps, and $P^-$ is a path
from $h_1$ to $h_0$ consisting only of down steps and loop steps. We
furthermore restrict ourselves to paths with the weight($P^+$)$\times$
weight($P^-$) f\/ixed to be $\by \mathbf w^\bn = w_1^{n_1}\cdots
w_{k}^{n_{k}}$ for some $\bn$. Here, $\by$ is the product of the
weights of the down steps in $P^-$, and $\bw^\bn$ is the total weight
coming from the loops in the path. Let $f$ be the weight of the
remaining f\/ixed portions of the path.

Without loss of generality, we can take $h_0=0$ and $h_1=k+1$.

For each such path we can decompose $n_i=n_i^++n_i^-$, where $n_i^+$
is the number of times the loop with weight $w_i$ is traversed in
$P^+$, and $n_i^-$ in $P^-$. Paths which arise from paths on
$\widetilde{H}_k$ must have $n_i^-\geq 1$ (for all $i$) by
def\/inition. We claim that on $\w H'_k$, paths with $n_i^-=0$ are
cancelled by paths which pass through the new ascending edges.

The key observation is that a path which has an up step going through
the ascending oriented edge $i-1\to i+a+1$ ($a\geq 0$) on $\w H'_k$ has
$n_i^+=n_{i+1}^+=\cdots=n_{i+a}=0$.

Then the total contribution of the paths in $\mathcal P$ to the
partition function is in fact
\begin{gather*}
f \by \sum_{n_i^++n_i^-=n_i\atop n_i^->0} \mathbf w^\bn  =
f \by \sum_{n_i^++n_i^-=n_i\atop n_i^+>0} \mathbf w^\bn \\
\hphantom{f \by \sum_{n_i^++n_i^-=n_i\atop n_i^->0} \mathbf w^\bn}{}
 =  f \by \sum_{n_i^++n_i^-=n_i}\mathbf w^\bn
-
 f \by \sum_{j=1}^k\sum_{n_i^++n_i^-=n_i\atop n_j^+=0}\mathbf w^\bn  \\
\hphantom{f \by \sum_{n_i^++n_i^-=n_i\atop n_i^->0} \mathbf w^\bn=}{}
+ f \by \sum_{j_1<j_2} \sum_{n_{j_1}^{+}=n_{j_2}^+=0}\mathbf w^\bn-
 f \by \sum_{j_1<j_2<j_3} \sum_{n_{j_1}^{+}=n_{j_2}^+=n_{j_3}^+=0}\mathbf
 w^\bn + \cdots\\
\hphantom{f \by \sum_{n_i^++n_i^-=n_i\atop n_i^->0} \mathbf w^\bn}{}
=  f \by \sum_{a=0}^k(-1)^a \sum_{j_1<\cdots<j_a}\sum_{n_{j_1}^+=\cdots
  =n^+_{j_a}=0} \mathbf w^{\bn}.
\end{gather*}
So the alternating sum has the ef\/fect of subtracting the terms with any
$n_i^-=0$.

A path on the graph $\w H'_k$ with $a$ spine vertices skipped (by
traversing ascending edges of length $>1$) comes with
a total sign $(-1)^x$ where $x=\sum(\hbox{length of the ascending
segments}-1)=a$. That is also the sign of the term with $n_{j_1}^+= \cdots
= n_{j_a}^+=0$ in the summation.

Finally, we note that any path can be decomposed into pairs of
ascending and descending segments as above, and the proof can be
applied iteratively to any path.
\end{proof}

\subsection{Def\/inition of compactif\/ied graphs}
On a graph $\Gamma_\bm$, we call a skeleton edge {\em horizontal} if it
connects (a) vertices~$i$ and $i'$ for some~$i$, (b)~vertices $0$ and
$1$, or (c)~the top vertex and the one below it. We call an edge
{\em vertical} otherwise.

\begin{definition}
The compactif\/ied graph $\Gamma_{\bm}'$ with $r+1$ vertices is obtained
from the graph $\Gamma_\bm$ via the following compactif\/ication
procedure:
\begin{enumerate}\itemsep=0pt
\item Introduce an order on the vertices of $\Gamma_\bm$, so that
  $i<i+1$ and $i<i'<i+1$. Number them from $1$ to $2r+2$ accordingly.
\item Identify vertices $2j-1$, $2j$ ($j=1,\dots ,r+1$), and rename the
  resulting vertex $j$. Double edges connecting $(2j-1,2j)$ are
  replaced by a loop at $j$ with weight which is the product of the
  weights on the two edges. All other edges and their weights are unchanged.
\item All maximal subgraphs of the form $\w H_k$ of $\Gamma_{\bm}$,
  consisting of vertical edges only, are replaced by compactif\/ied
  weighted graphs of the form $\w H_k'$, as in Lemma~\ref{inclexcl} (with the
  obvious shift in labels).
\end{enumerate}
\end{definition}
\begin{figure}[t]
\centerline{\includegraphics{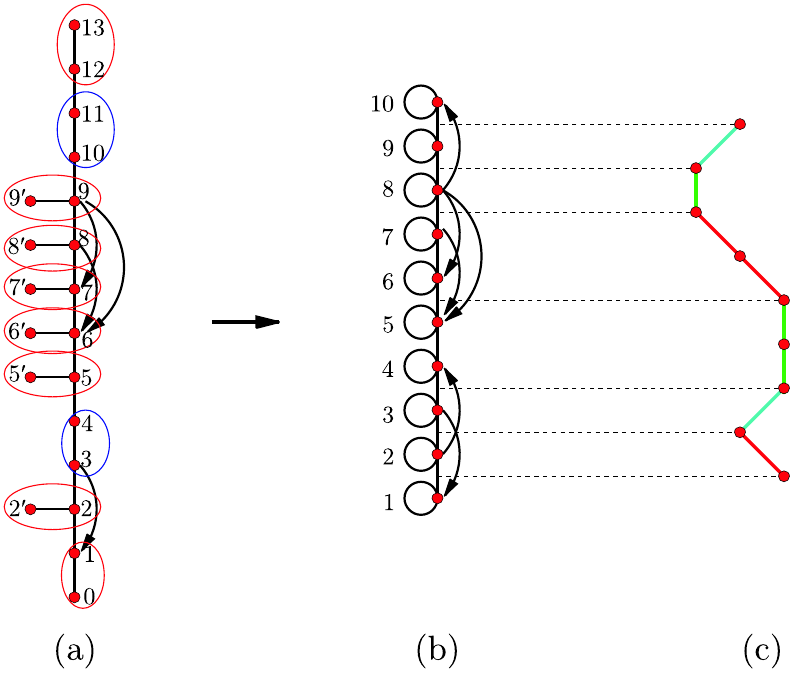}}
\caption{(a) The graph $\Gamma_\bm$ of Fig.~\ref{fig:exglue},
and (b) the compacted graph $\Gamma_\bm'$. We have circled on
$\Gamma_\bm$ the pairs of vertices to be identif\/ied in the
compactif\/ication procedure (in red for horizontal pairs, in blue for
vertical pairs). We also represent in (c) the Motzkin path $\bm$.}
\label{fig:excompact}
\end{figure}

\begin{example}
For illustration, the identif\/ication of edges in the case of the graph
of Fig.~\ref{fig:exglue}~(c) are: $0\sim 1$, $2\sim 2'$, $3\sim 4$,
$5\sim 5'$, $6\sim 6'$, $7\sim 7'$, $8\sim 8'$, $9\sim 9'$, $10\sim
11$ and $12\sim 13$. See Fig.~\ref{fig:excompact}. There are two maximal
subgraphs of the form $\w H_1$ are the vertices $2$, $3$, $4$, $5$ and
$9$, $10$, $11$, $12$.  Each pair now corresponds to a vertex~$i$ in
$\Gamma_\bm'$, $i=1,2,\dots ,10$, which receives a loop $i\to i$ from the
identif\/ication.
\end{example}

\begin{figure}[t]
\centerline{\includegraphics{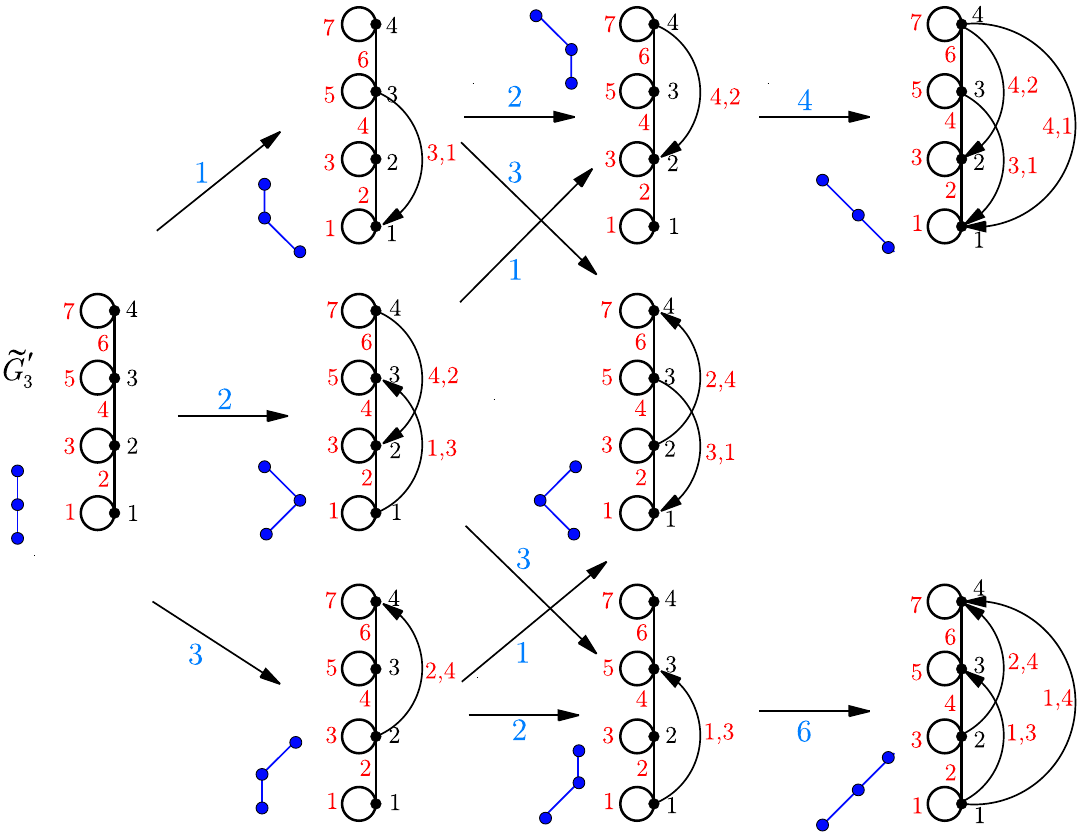}}
\caption{The Motzkin paths $\bm$ of the fundamental domain ${\mathcal M}_3$
and the associated compact graphs $\Gamma_\bm'$,
with their vertex and edge labels. Note the up-step weights:
$y_{1,3}=y_{2,4}=-1$ and $y_{1,4}=1$.
Mutations are indicated by arrows.}\label{fig:compmotzfour}
\end{figure}

Fig.~\ref{fig:compmotzfour} shows the set of compactif\/ied
graphs corresponding to Fig.~\ref{fig:pathmotzfour} for the case~$A_3$.

The resulting weighted graph $\Gamma_\bm'$ has $r+1$ vertices labelled
$1,2,\dots ,r+1$, and hence is associated with a transfer matrix $T_\bm'$
of size $r+1 \times r+1$.

\subsection{An alternative construction}\label{section6.4}
An alternative description of the compactif\/ied graphs $\Gamma'_\bm$ is
the following.

We start from the graph $\Gamma_{\bm_0}'\equiv{\widetilde G}_r'$ of
Fig.~\ref{fig:dualgrprime}. The loop at vertex $i$ has weight $t
y_{2i-1}$ and the edge $i+1\to i$ has weight $t y_{2i}$, where
$y_j=y_j(\bm_0)$ are as in equations \eqref{oddy}, \eqref{eveny}.

Decompose $\bm$ into maximal segments of the form:{\samepage
\begin{enumerate}\itemsep=0pt
\item[$1)$] descending segments, $D_{\al,i}=((\al,m),(\al+1,m-1),\dots ,(\al+i-1,m-i+1))$;
\item[$2)$] ascending segments,
  $A_{\al,i}=((\al,m),(\al+1,m+1),\dots ,(\al+i-1,m+i-1))$;
\item[$3)$] f\/lat segments $((\al,m),(\al+1,m),\dots ,(\al+k-1,m))$.
\end{enumerate}
Here, $i\geq 2$ and $k\geq 1$.}

\begin{definition}\label{gammapp}
The graph $\Gamma_\bm''$ is the graph obtained from $\Gamma_{\bm_0}'=\w G_r'$
via the following steps:
\begin{enumerate}\itemsep=0pt
\item For each descending sequence $D_{\al,i}$ we add descending edges
$\al+p\to \al+q$  ($0< q+1 < p \leq i$) to $\Gamma_{\bm_0}'$, with
weights $t y_{\al+p,\al+q}(\bm)$, where
\begin{equation}\label{redunwei}
y_{\al+p,\al+q}(\bm)={\prod\limits_{j=\al+q}^{\al+p-1}y_{2j}(\bm)
\over \prod\limits_{j=\al+q+1}^{\al+p-1} y_{2j-1}(\bm) }.
\end{equation}
\item
For each ascending sequence $A_{\al,i}$, we add ascending
edges $\al+q \to \al+p$ ($0<q+1<p\leq i$) with
weights
$y_{\al+q,\al+p}(\bm)=(-1)^{p-q-1}$.
\end{enumerate}
\end{definition}

\begin{lemma}\label{ident}
The weighted graph $\Gamma''_\bm$ is identical to the weighted graph
$\Gamma'_\bm$.
\end{lemma}
This is just the result of the def\/inition of $\Gamma_\bm$ using the
decomposition of $\bm$, as in Fig.~\ref{fig:exglue}. Maximal subgraphs the
form $\w H_k$ correspond to the maximal ascending segments of the
Motzkin path. All other segments correspond to subgraphs with
horizontal edges.

\subsection{Equality of generating functions}
To summarize, the compactif\/ied graph $\Gamma'_\bm$ is such that
\begin{theorem}\label{eqgenfns}
\[
(1-T_{\bm})^{-1}_{1,1}=(1-T_{\bm}')^{-1}_{1,1}.
\]
\end{theorem}
In other words: the partition function for weighted paths from vertex
$1$ to vertex $1$ in $\Gamma_{\bm}$ is identical to that for weighted
paths from vertex $1$ to vertex $1$ in the compact graph
$\Gamma_\bm'$.

\section{Totally positive matrices and compactif\/ied transfer
  matrices}\label{totalpos}

We now establish the connection between the transfer matrices $T'_\bm$
for paths on $\Gamma'_\bm$ and the totally positive matrices of
\cite{FZposit} corresponding to double Bruhat cells for pairs of
Coxeter elements.

We may express the compact transfer matrices $T_\bm'$ of the previous
section in terms of the elementary matrices $f_i$, $e_i$, $d_i$ for
$GL_{r+1}$, def\/ined as follows. Let $E_{ij}$ denote the standard elementary
matrix of size $(r+1)\times(r+1)$, with entries
$(E_{i,j})_{k,\ell}=\delta_{k,i}\delta_{\ell,j}$.
\begin{definition}The elementary matrices $\{e_i, f_i, d_i\}$ are def\/ined by
\begin{gather}
f_i = I +\lambda_i E_{i+1,i}, \qquad e_i=I+\nu_i E_{i,i+1},
\qquad i\in\{1,\dots ,r\},\nonumber\\
  d_i=I + (\mu_i-1)E_{i,i},\qquad i\in \{1,\dots ,r+1\}.\label{elementary}
\end{gather}
for some real parameters $\lambda_i$, $\mu_i$, $\nu_i$.
\end{definition}
In \cite{FZposit}, Fomin and Zelevinsky introduced a parametrization
of totally positive matrices as products of the form $\prod\limits_{i\in I}
f_i \prod\limits_{i=1}^{r+1} d_i \prod\limits_{j\in J}e_j$ for $I$, $J$ two suitable
sets of indices, and  $\lambda_i$, $\mu_i$, $\nu_i$ some positive parameters.
This expression allowed to rephrase total positivity
in terms of networks.  Here we interpret our compact transfer matrices
in terms of some of these products.

Recall that each Motzkin path can be decomposed into descending,
ascending and f\/lat pieces, as in Section~\ref{section6.4}. We introduce the
increasing sequence of integers $(a_1,\dots ,a_{2k})$, such that the
$j$th ascending piece of $\bm$, $A_{{\al_j},i_j}$ of $\bm$ starts at
$m_{\al_j}=a_{2j-1}$ and ends at $m_{{\al_j}+i_j-1}=a_{2j}$. Similarly
for the sequence of increasing integers $(b_1,\dots ,b_{2p})$, which mark
the starting and ending points of the descending sequences
$D_{\al_j,i}$.

For $i<j$, let $\omega[i,j]$ denote the permutation which reverses the
order of all consecutive elements between $i$ and $j$ in a given
sequence. That is,
$\omega[i,j]=(i,j)(i+1,j-1)(i+2,j-2)\cdots $.
For example,  $\omega[i,j]$: $(j,j-1,\dots ,i)\mapsto (i,i+1,\dots ,j)$,
\begin{gather}
\sigma_\bm  =  \left(\prod_{i=1}^k\omega[a_{2i-1},a_{2i}]\right) \circ
(r,r-1,\dots ,1),\nonumber \\
\tau_{\bm}  =   \left( \prod_{i=1}^p\omega[b_{2i-1},b_{2i}]\right) \circ
(r,r-1,\dots ,1).\label{rearrangement}
\end{gather}

\begin{example}\label{examotz}
For the Motzkin path $\bm=(2,1,2,2,2,1,0,0,1)$ of
Fig.~\ref{fig:exglue}, we have the ascending segments $[2,3]$ and
$[8,9]$, while the descending segments are $[1,2]$ and $[5,7]$. The
rearranged sequences read $\sigma_\bm=(8,9,7,6,5,4,2,3,1)$ and
$\tau_\bm=(9,8,5,6,7,4,3,1,2)$.
\end{example}

Note that the sequences $\sigma_\bm$ and $\tau_\bm$ consist of
increasing and decreasing subsequences of consecutive integers, and
that these subsequences and their order are unique.

One can def\/ine the decomposition of the transfer matrix $T_\bm'$ into
a strictly lower-triangular part $N_\bm$ and an upper triangular part
$B_\bm$, so that
\[
T_{\bm}' =N_{\bm} +B_\bm.
\]
\begin{lemma}\label{decomposition}
The matrices $N_\bm$ and $B_\bm$ can be expressed as
\begin{gather}
N_\bm = I-(f_{i_1}f_{i_2}\cdots f_{i_r})^{-1},
\label{lower}\\
B_\bm = t (d_1d_2 \cdots d_{r+1})\,
(e_{j_1}e_{j_2}\cdots e_{j_r}),\label{upper}
\end{gather}
where
the parameters in equation \eqref{elementary} are
\begin{equation}\label{munu}
\lambda_i=1,\qquad \mu_i=y_{2i-1},\qquad \nu_i={y_{2i}\over y_{2i-1}}.
\end{equation}
\end{lemma}
\begin{proof}
We give a pictorial proof. It is possible to describe multiplication
by an elementary matrix as the addition of an arrow to a graph.  In
our context, $N_\bm$ encodes the ascending arrows in the
graph~$\Gamma_\bm'$, and~$B_\bm$ the descending arrows.

First, consider the product $f_{i_r}^{-1}\cdots f_{i_1}^{-1}$ in
$N_\bm$ of equation~\eqref{lower}, where $(f_i)^{-1}=I-E_{i+1,i}$.

The sequence $(i_r,i_{r-1},\dots ,i_1)=\big(\prod_i
\omega[a_{2i-1},a_{2i}]\big)\circ(1,\dots ,r)$, which is simply $\sigma_\bm$
written in reverse order, consists of alternating increasing and
decreasing sequences of consecutive integers. The products of matrices
corresponding to increasing subsequences are
\begin{gather*}
P^+_{j} =  f^{-1}_{a_{2j}+1} f^{-1}_{a_{2j}+2} \cdots
f^{-1}_{a_{2j+1}-1}
= I - \sum_{i=1}^{a_{2j+1}-a_{2j}-1}
E_{a_{2j}+i+1,a_{2j}+i}
\end{gather*}
and products corresponding to decreasing sequences are
\begin{gather*}
P^-_j =  f^{-1}_{a_{2j}} f^{-1}_{a_{2j}-1} \cdots f^{-1}_{a_{2j-1}}
=I+ \sum_{i=1}^{a_{2j}-a_{2j-1}} \sum_{k=0}^{i-1}
(-1)^{i+k}E_{a_{2j-1}+i,a_{2j-1}+k}.
\end{gather*}

We start with the graph corresponding to the identity matrix, which is
the transfer matrix of the graph consisting of $r+1$ disconnected
vertices labelled $1,2,\dots ,r+1$, each with a loop of weight $1$.
Multiplying on the left by $f_j^{-1}$ creates an ascending edge $j\to
j+1$ with weight $-1$. More generally, left multiplication by $P^+_j$
creates a succession of ascending edges $a_{2j}+1\to a_{2j}+2$,
$a_{2j}+2\to a_{2j}+3$, \dots,  $a_{2j+1}-1\to a_{2j+1}$.  Left
multiplication by $P^-_j$ creates a web of ascending edges
$a_{2j-1}+k\to a_{2j-1}+i$, $0\leq k\leq i-1\leq a_{2j}-a_{2j-1}-1$,
with alternating weights $(-1)^{i+k}$.  We illustrate the resulting
actions on the graphs in Fig.~\ref{fig:incdec}.

\begin{figure}[t]
\centerline{\includegraphics{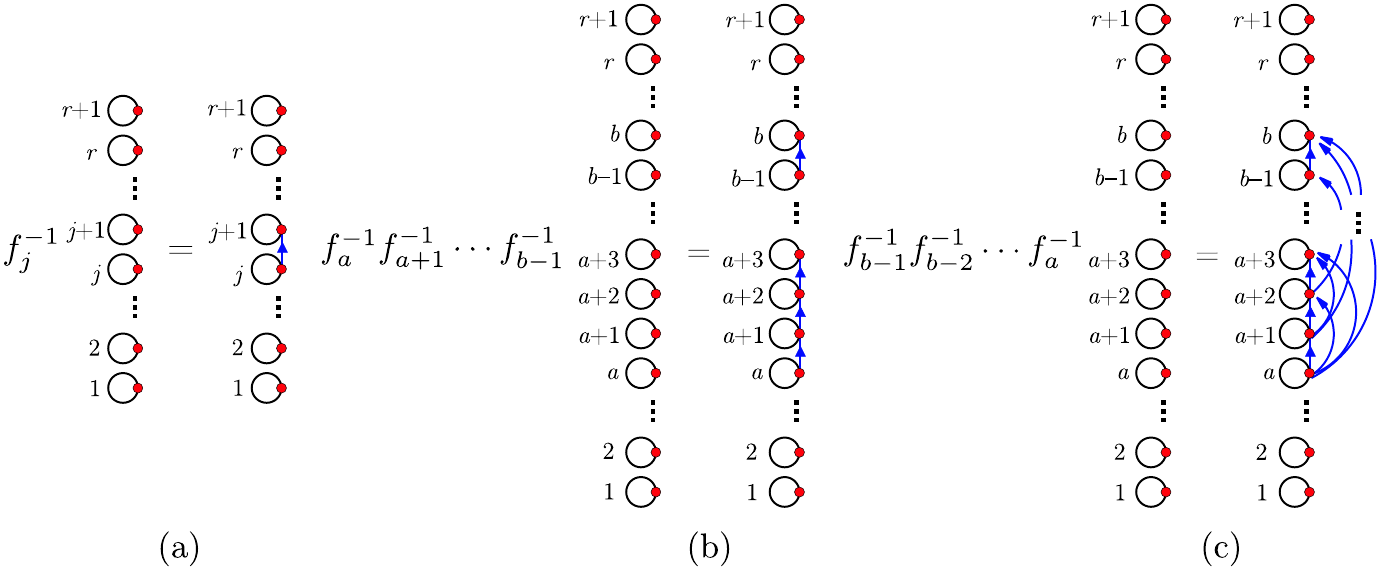}}
\caption{Pictorial representation of the action of (a) $f_j^{-1}$,
(b) an increasing product $f_a^{-1}f_{a+1}^{-1} \cdots f_{b-1}^{-1}$,
and (c) a decreasing product $f_{b-1}^{-1}f_{b-2}^{-1}\cdots f_a^{-1}$
for some $a<b$.
These correspond to adding (a) an ascending edge $j\to j+1$ with weight~$-1$,
(b) ascending edges $a+i\to a+i+1$, $i=0,1,\dots ,b-a-1$ with weights~$-1$,
and (c) a web of ascending edges $a+i\to a+k$, $0\leq i\leq k-1\leq b-a-1$
with weights $(-1)^{i+k}$.}
\label{fig:incdec}
\end{figure}

Recall that the segments $[a_{2i-1},a_{2i}]$ correspond to the
ascending segments of $\bm$, themselves associated to the vertical
chain-like pieces of $\Gamma_\bm$ (see Fig.~\ref{fig:exglue}). Recall that in
the identif\/ication procedure leading to $\Gamma_\bm'$ (Def\/inition
\ref{gammapp} and Lemma \ref{ident}),
we showed that such chains must receive a web of ascending edges
with alternating weights $\pm 1$ (see Fig.~\ref{fig:loopsit}~(c)),
while all the vertices are connected via ascending edges $i\to i+1$
with weight $1$.

Finally, comparing this with the graph associated to
$(f_{i_1}f_{i_2}\cdots f_{i_k})^{-1}$ as described above, we f\/ind that
the contribution of ascending edges to $T_\bm'$ (or equivalently,
$\Gamma_\bm'$) is identical to $I-(f_{i_1}f_{i_2}\cdots
f_{i_k})^{-1}$.

The proof of equation \eqref{upper} is similar, but now concerns
the descending edges and the loops of the graph $\Gamma_{\bm}'$.

The product $t d_1d_2\cdots d_{r+1}$ is the transfer matrix of
a chain of $r+1$ disconnected vertices $i=1,2,\dots ,r+1$, each with a
loop with weight $ty_{2i-1}$.
Multiplication on the right by $e_i$ creates a descending edge
$i+1\to i$, with weight $ty_{2i-1} {y_{2i}\over y_{2i-1}} =t
y_{2i}$.

Again we divide the sequence $\tau_\bm$ into increasing
subsequences of consecutive integers,
$(b_{2j-1}$, $b_{2j-1}+1,\dots ,b_{2j})$ and decreasing subsequences
$(b_{2j+1}-1,\dots ,b_{2j}+1)$. Therefore the product
$e_{j_1}e_{j_2}\cdots e_{j_r}$ consists of ``ascending'' factors
\[
Q^+_j=e_{b_{2j-1}}e_{b_{2j}}\cdots e_{b_{2j}}=
I+\sum_{i=1}^{b_{2j}-b_{2j-1}} \sum_{k=0}^{i-1}
y_{b_{2j-1}+i,b_{2j-1}+k}E_{b_{2j-1}+k,b_{2j-1}+i}
\]
and ``descending'' factors
\[
Q^-_j=e_{b_{2j+1}-1} e_{b_{2j+1}-2}\cdots e_{b_{2j}+1}=
I+\sum_{i=1}^{b_{2j+1}-b_{2j}-1}
\frac{y_{2b_{2j}+2i}}
{y_{2b_{2j}+2i-1}}\, E_{b_{2j}+i,b_{2j}+i+1},
\]
where $y_{b+i, b+j}$ are the weights of equation \eqref{redunwei}.

Recall that the segments
$[b_{2i-1},b_{2i}]$ correspond to the descending segments of $\bm$,
which correspond to networks of descending edges on $\Gamma_\bm$ with
weights \eqref{redunwei} (see Fig.~\ref{fig:exglue}).  In our
construction \ref{gammapp} of $\Gamma_\bm'$, these descending
edges have remained unchanged, while each vertex $i$ received a loop
with weight $ty_{2i-1}$, $i=1,2,\dots ,r+1$. This is
nothing but the graph associated to $t (d_1d_2\cdots
d_{r+1})(e_{j_1}e_{j_2}\cdots e_{j_r})$, which therefore encodes the
contribution of loops and descending edges to $T_\bm'$, and
equation~\eqref{upper} follows.
\end{proof}

We can now make a direct connection with totally positive matrices encoding
the networks associated to the Coxeter double Bruhat cells considered
in \cite{GSV}.  Each Motzkin path $\bm\in {\mathcal M}_r$ corresponds
to such an element, the $(r+1)\times (r+1)$ matrix $P_\bm$:
\begin{definition} Given a Motzkin path $\bm\in \mathcal M_r$, def\/ine
\begin{equation*}
P_{\bm}= (f_{i_1}f_{i_2}\cdots f_{i_r})\, (d_1 d_2 \cdots d_{r+1})\,
(e_{j_1}e_{j_2}\cdots e_{j_r}),
\end{equation*}
where $(i_1i_2\cdots i_r)=\sigma_\bm$ and $(j_1 j_2\cdots
j_r)=\tau_\bm$ are the two sequences of \eqref{rearrangement}.
The parameters $\lambda_j$, $\mu_j$, $\nu_j$ are as in equation \eqref{munu}.
\end{definition}

As a consequence of equations \eqref{lower}, \eqref{upper}, we have
\begin{theorem}
\begin{equation*}
(I-T_\bm' )^{-1}=(I-t P_\bm)^{-1}\, (f_{i_1}f_{i_2}\cdots f_{i_r}),
\end{equation*}
which allows to rewrite the generating function of Theorem {\rm \ref{eqgenfns}}
\begin{equation*}
\big((I-T_\bm' )^{-1}\big)_{1,1}=\left\{ \begin{matrix}
 \big((I-t P_\bm)^{-1}\big)_{1,1}+\big((I-t P_\bm)^{-1}\big)_{1,2}  & \ \ {\rm if} \ \ a_1>1,\vspace{1mm} \\
\displaystyle \sum_{a=a_1}^{a_2+1}  \big((I-t P_\bm)^{-1}\big)_{1,a} &
\ \ {\rm if} \ \ a_1=1. \end{matrix}\right.
\end{equation*}
\end{theorem}
This yields an interpretation of the solution $R_{1,n+m_1+1}$
to the $A_r$ $Q$-system with initial data~$\bx_\bm$
in terms of the network associated to the totally
positive matrix $P_\bm$, for all $n\geq 0$.

\begin{example}
For the fundamental Motzkin path $\bm=\bm_0$
with $m_\al=0$ for all $\al$, we
have $\sigma_{\bm_0}=\tau_{\bm_0}=(r,r-1,\dots ,1)$, and therefore
$P_{\bm_0}=(f_r f_{r-1} \cdots f_1)(d_1 d_2 \cdots d_{r+1})(e_r e_{r-1}\cdots e_1)$. Note that the matrix $F=f_r f_{r-1} \cdots f_1$
has entries $F_{i,j}=1$ if $i\geq j$, and $0$ otherwise. One
can check directly that $I-t P_{\bm_0}=F(I-T_{\bm_0}')$,
with $T_{\bm_0}'$ given by equation \eqref{transmatcompact}.
\end{example}

An equivalent formulation uses the explicit decomposition $P_\bm=FDE$,
where $F=f_{i_1}\cdots f_{i_r}$, $D=d_1\cdots d_{r+1}$ and
$E=e_{j_1}\cdots e_{j_r}$:
\begin{equation*}
(I-t P_\bm)^{-1} F=\sum_{n=0}^\infty t^n (FDE)^n F=
\sum_{n=0}^\infty t^n F(DEF)^n =F(I-t P_\bm')^{-1},
\end{equation*}
where $P'_\bm=DEF=F^{-1}P_\bm F$, and
the fact that $F$ is a lower uni-triangular matrix, which implies:
$\big((I-t T'_\bm)^{-1}\big)_{1,1}=\big((I-t P_\bm)^{-1} F\big)_{1,1}=
\big((I-t P'_\bm)^{-1}\big)_{1,1}$.

\begin{figure}[t]
\centerline{\includegraphics[scale=1.05]{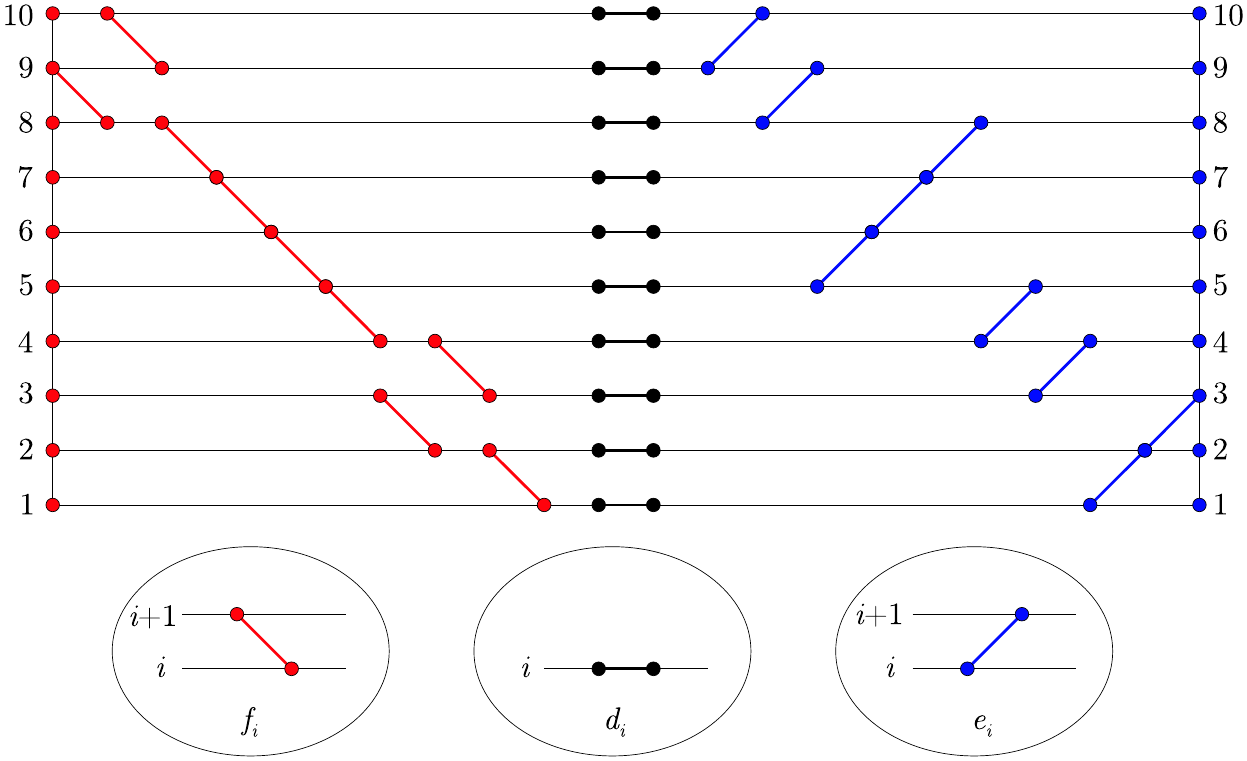}}
\caption{The network corresponding to the totally positive matrix
$P_{\bm}$ for the Motzkin path $\bm$ of Example~\ref{examotz},
depicted in Fig.~\ref{fig:excompact} (c). We have indicated in medallions the
three network representations (in red, black, blue)
for the elementary matrices $f_i$, $d_i$ and $e_i$.
The network for $P_\bm'$ corresponds to the same picture, but with the
$f$ part (red descending elements) to the right of the $e$ part (ascending blue elements).}
\label{fig:exnetwork}
\end{figure}

The matrix $P'_\bm$ is another way to write a totally positive matrix,
and the network graph corresponding to it has a slightly modif\/ied form
from that of $P_\bm$. Both of these correspond to electrical networks
\cite{FZposit}.
For illustration, we represent in Fig.~\ref{fig:exnetwork}
the network corresponding to the matrix $P_\bm$ for
the Motzkin path $\bm$ of Example \ref{examotz},
represented in Fig.~\ref{fig:excompact} (c). The medallions show the
three elementary circuit representations for the three types of elementary
matrices $f_i$, $d_i$, $e_i$, each receiving the associated weight. The Lindstr\"om
lemma \cite{LGV1} of network
theory states that the minor $|P|_{r_1,\dots ,r_k}^{c_1,\dots ,c_k}$
of the matrix $P$ of the network,
corresponding to  a specif\/ic choice rows $r_1,\dots ,r_k$ and columns $c_1,\dots ,c_k$,
is the partition function of $k$ non-intersecting (vertex-disjoint) paths
starting at points $r_1,\dots ,r_k$ and ending at points $c_1,\dots ,c_k$, and with steps
taken only on horizontal lines or along $f$, $d$ or $e$ type elements. Here we
have only considered circuits with entry and exit point $1$, after possibly several iterations
of the same network (each receiving the weight $t$),
and whose generating function is precisely the resolvent
$\big((I-t P)^{-1}\big)_{1,1}$.

\section{Conclusion}\label{conc}

\begin{figure}
\centerline{\includegraphics{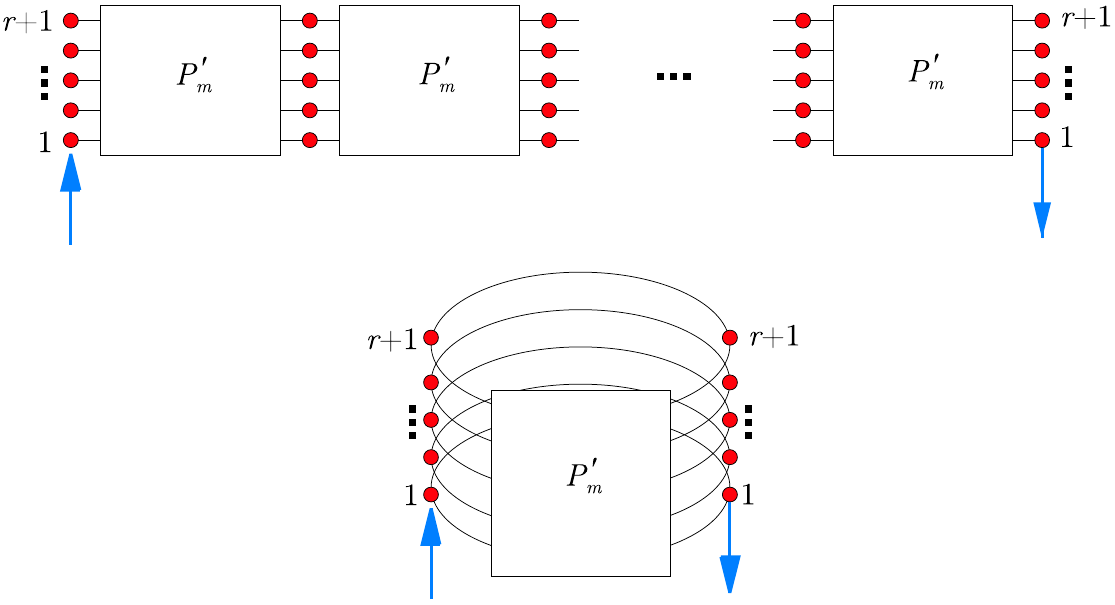}}
\caption{The concatenation of $n$ copies of the network coded by $P_\bm'$.
The quantity $(P_\bm')^n_{1,1}$ is the partition function of electrical wires starting
and ending at the two indicated arrows. We also represent below a cylinder formulation
\`a la \cite{GSV}: the wire must wind $n$ times around the cylinder before exiting.}
\label{fig:concanet}
\end{figure}

In this paper, we have made the contact between our earlier study of the
solutions of the $A_r$ $Q$-system, expressed in terms of initial data
coded via Motzkin paths, and the totally positive matrices for Coxeter
double Bruhat cells. We showed in particular how the relevant
pairs of Coxeter elements were encoded in the Motzkin paths as well.

One would expect the total positivity of the transfer matrices $P_\bm$
or $P_\bm'$ to be directly related to the proof of the positivity
conjecture in the case of the $A_r$ $Q$-system.  Our proof presented
in \cite {DFK08a} relies on the path formulation of $R_{1,n}$ and on
the formulation of $R_{\al,n}$ as the partition function of families
of strongly non-intersecting paths. The total positivity of the
compactif\/ied formulation should provide an alternative proof, using
networks rather than paths.

The precise connection between paths on graphs and networks, as
illustrated in Section \ref{totalpos} above is subtle. Indeed, the
identity between resolvents implies that the partition function for
weighted paths from $1$ to $1$ on $\Gamma_\bm$ with $n$ descents,
$(T_\bm^n)_{1,1}$, is identical to the generating function for
circuits on a network made of $n$ identical concatenated networks,
each corresponding to the totally positive matrix $P'_\bm$, from
connector $1$ to connector $1$ (see the top of Fig.~\ref{fig:concanet}).
In \cite{GSV}, this concatenation is realized by putting the network on
a cylinder and allowing for the circuit to wind $n$ times around it
before exiting (see the bottom of Fig.~\ref{fig:concanet}). Note
that we could also work with $P_\bm$ instead, as it is related to
$P_\bm'$ via cyclic symmetry.

More generally, it should be possible to relate our non-intersecting
path families to networks with multiple entries and exits, as in the
setting of the Lindstr\"om lemma.

Another question concerns the cluster algebra attached to the $A_r$
$Q$-system.  As stressed in~\cite{DFK08a}, we have only considered a
subset of the clusters which arise in the full $Q$-system cluster
algebra, namely those which consist of solutions of the $Q$-system.
There are other cluster mutations, however, which are not recursion
relations of the form \eqref{qsys}.  One may ask about the other
cluster variables in the algebra. The positivity conjecture should
hold for them as well. Preliminary investigations show that the
corresponding mutations can still be understood in terms of (f\/inite)
continued fraction rearrangements, hence we expect them to also have a~network counterpart. These clearly can no longer correspond to Coxeter
double Bruhat cells, as those are exhausted by the solutions of the
$Q$-system.

Finally, the connection to total positivity should be generalizable to
the case of other simple Lie algebras as well. Indeed, on the one hand the
$Q$-systems based on other Lie algebras also have cluster algebra
formulations \cite{DFK08}, while on the other hand the notion of total
positivity has been extended to arbitrary Lie groups \cite{FZpositII}.
We have evidence that hard particle and path interpretations exist for
all $Q$-systems, and it would be interesting to investigate their
relation to the corresponding generalized networks. The integrability
of these systems is presumably related to that of the Coxeter--Toda
systems of \cite{RESH}. This will be the subject of forthcoming work.

\subsection*{Acknowledgements}

We thank M.~Gekhtman, S.~Fomin, A.~Postnikov, N.~Reshetikhin
and A.~Vainshtein for useful discussions.
RK's research is funded in part by NSF grant DMS-0802511.
RK thanks CEA/Saclay
IPhT for their hospitality. PDF's research is partly supported by the
European network grant ENIGMA and the ANR grants GIMP and GranMa.
PDF thanks the department of Mathematics of the University
of Illinois at Urbana-Champaign for hospitality and support, and the
department of Mathematics of the University of California Berkeley
for hospitality.

\addcontentsline{toc}{section}{References}
\LastPageEnding

\end{document}